\documentclass[11pt]{article}
\usepackage[T1]{fontenc}
\usepackage[utf8]{inputenc}
\usepackage{lmodern}
\usepackage[a4paper,margin=28mm]{geometry}
\usepackage{authblk}

\usepackage{amsmath,amssymb,amsthm,mathtools}
\usepackage{microtype,booktabs,array,enumitem}
\usepackage{tikz}
\usepackage[numbers,sort&compress]{natbib}
\usepackage[hidelinks]{hyperref}
\hypersetup{pdftitle={Positive diagonal congruence and weight transformations of Ihara zeta functions},
 pdfsubject={English mathematical research-paper working draft},
 pdfkeywords={Ihara zeta, diagonal congruence, cavity equations, nonnegative fixed points}}
\usepackage[nameinlink,noabbrev]{cleveref}
\usepackage{fancyhdr}
\setlist[enumerate]{itemsep=3pt,topsep=5pt}
\numberwithin{equation}{section}
\theoremstyle{plain}
\newtheorem{theorem}{Theorem}[section]
\newtheorem{proposition}[theorem]{Proposition}
\newtheorem{lemma}[theorem]{Lemma}
\newtheorem{corollary}[theorem]{Corollary}
\theoremstyle{definition}
\newtheorem{definition}[theorem]{Definition}
\newtheorem{example}[theorem]{Example}
\newtheorem{algorithm}{Algorithm}[section]
\crefname{algorithm}{Algorithm}{Algorithms}
\Crefname{algorithm}{Algorithm}{Algorithms}
\theoremstyle{remark}
\newtheorem{remark}[theorem]{Remark}
\newcommand{\RR}{\mathbb R}
\newcommand{\CC}{\mathbb C}
\newcommand{\uu}{\boldsymbol u}
\newcommand{\vv}{\boldsymbol v}
\newcommand{\one}{\boldsymbol 1}
\newcommand{\dd}{\boldsymbol\delta}
\newcommand{\cav}{\boldsymbol c}
\newcommand{\zz}{\boldsymbol z}
\newcommand{\svec}{\boldsymbol s}
\newcommand{\arcs}{E_D}
\DeclareMathOperator{\diag}{diag}
\DeclareMathOperator{\Spec}{Spec}
\DeclareMathOperator{\Tr}{Tr}

\DeclareMathOperator{\im}{im}

\title{Positive diagonal congruence and weight transformations\ of Ihara zeta functions}
\author[1]{Yuma Furuta\thanks{furuta.yuma@phys.kyushu-u.ac.jp}}
\author[2]{So Matsuura\thanks{s.matsu@keio.jp}}
\affil[1]{Institute for Advanced Study, Kyushu University, Fukuoka 812-8581, Japan}
\affil[2]{Hiyoshi Departments of Physics, and Research and Education Center for Natural Sciences, Keio University, Yokohama, Kanagawa 223-8521, Japan}
\date{November 2026}
\begin{document}

\maketitle
\begin{abstract}
We construct identities relating the values of the Ihara zeta function at two sets of independent edge weights in $(0,1)$ on the same finite graph. Positive diagonal congruences between weighted vertex matrices in the Ihara--Bass formula are parameterized by positive solutions of the Gaussian cavity equations with positive associated vertex coefficients. When the non-backtracking edge adjacency matrix is irreducible, a unique such cavity solution exceeding an explicit reference solution in every component exists exactly when that edge adjacency matrix, weighted by the squared input edge weights, has spectral radius greater than one. This cavity solution determines a real-analytic weight transformation and an explicit prefactor relating the input and output zeta functions. The non-backtracking edge adjacency matrix weighted by the squared output edge weights has spectral radius below one. Positive definiteness of the weighted vertex matrix permits only the identity weight transformation. On a chain of $k$ labelled copies of $K_4$ joined by sufficiently weak edges, we prove that there are at least $2^k$ distinct positive cavity solutions with positive vertex coefficients. 
For uniform weights, the resulting functional equations recover the regular and biregular reciprocal relations.
\end{abstract}
\noindent\textbf{Keywords.} Ihara zeta function; edge adjacency matrix; diagonal congruence; Gaussian cavity equations; nonnegative fixed points.

\smallskip
\noindent\textbf{MSC 2020.} Primary 05C50; Secondary 05C38, 15A18, 15B48, 92D30.

\smallskip
\noindent\textbf{Preprint.} KYUSHU-HET-388

\section{Introduction}\label{sec:introduction}
\subsection{The problem: changing weights without changing the graph}

Ihara introduced his zeta function in an arithmetic setting \cite{Ihara1966}. Serre observed that it admits a graph-theoretic interpretation as an Euler product over primitive reduced cycles \cite{Serre1980}. On a regular graph of degree at least two, the zeta function satisfies a classical functional equation relating two reciprocal values of its argument. We ask whether a change of independent edge weights on a fixed irregular graph can yield an analogous identity.

Sunada and Hashimoto developed the graph-theoretic theory \cite{Sunada1986,Hashimoto1989,Hashimoto1990}. The reciprocal zeta function is the determinant of $I-uW$, where $W$ is the non-backtracking matrix on oriented edges. The Ihara--Bass formula expresses this determinant in terms of a matrix on vertices and an explicit edge factor \cite{Hashimoto1989,Bass1992}. For other proofs see \cite{StarkTerras1996,FoataZeilberger1999,KotaniSunada2000}. Bartholdi extended the cycle count by introducing a parameter for backtracking \cite{Bartholdi1999}.

Our edge weights are scalar and take the same value on the two orientations of an edge. The multivariable formulation goes back to the edge zeta functions of Stark and Terras \cite{StarkTerras1996}. Weighted determinant and covering formulas were developed in \cite{MizunoSato2004,HortonStarkTerras2008}, and Artin-type factorization for graph coverings is treated in \cite{StarkTerras2000}. Morita gave a common framework for several weighted graph zeta functions \cite{Morita2020}. The corresponding vertex-determinant formulas were established for finite simple graphs \cite{IdeEtAl2021} and extended to finite digraphs with different inverse-arc conventions \cite{IshikawaMoritaSato2022,IshikawaMorita2026}. For an account of graph zeta functions see \cite{Terras2011}, and for a survey of their combinatorial expressions in Japanese see \cite{MoritaIshikawa2025Survey}.

For a $(t+1)$-regular graph with $t>0$ and adjacency matrix $A$, the ordinary Ihara--Bass formula reads
\[
 \zeta_G(u)^{-1}=(1-u^2)^{|E|-|V|}\det(I-uA+tu^2I).
\]
Replacing $u$ by $1/(tu)$ multiplies the matrix inside the determinant by $1/(tu^2)$. Together with the edge factor this gives the classical reciprocal relation \cite{Hashimoto1989,Terras2011}. On an irregular graph the coefficient of $u^2$ is the diagonal matrix of degrees minus one, so this scalar rescaling is no longer available in general.

We use \emph{functional equation} for a zeta identity induced by a prescribed transformation of the edge weights on the same graph. Both the transformed weights and the prefactor are determined by local equations on oriented edges, independently of the values of the zeta function. The classical reciprocal transformation is an involution. For the general transformation constructed below, the domain and the location of its image are part of the result.

\subsection{Notation}
\label{sec:notations}
Throughout this paper, 
$G=(V,E)$ is a finite connected simple undirected graph with at least two vertices. Write $n=|V|$, $m=|E|$, and let $\arcs$ contain the two orientations of every edge, which are called \emph{arcs}. For an arc $e$, let $o(e)$ and $t(e)$ be its origin and terminus, and let $\bar e$ be the reversed arc.
We write $d_v$ for the degree of a vertex $v\in V$.

We use \emph{path} for a finite arc sequence $P=(e_1,\ldots,e_\ell)$ satisfying $t(e_i)=o(e_{i+1})$ for $1\le i<\ell$. In particular, we allow repeated vertices and edges for a path. 
A \emph{simple path} is a path with no repeated vertices.
A \emph{cycle} is a nonempty path with $t(e_\ell)=o(e_1)$. 
A \emph{simple cycle} is a cycle of length at least three with no repeated vertices apart from its coinciding endpoints.
A cycle $C=(e_1,\ldots,e_\ell)$ is \emph{reduced} if $e_{i+1}\ne\bar e_i$ for $1\le i<\ell$ and $e_1\ne\bar e_\ell$, so it has neither backtracking nor a tail.
A cycle is \emph{primitive} if it cannot be written as $D^r$ for a shorter cycle $D$ and an integer $r\ge2$. 
A \emph{cycle graph} is a connected simple graph in which every vertex has degree two.

The non-backtracking edge adjacency matrix (the Hashimoto matrix) is
\begin{equation}\label{eq:Wintro}
 W_{ef}:= \mathbf1_{\{t(e)=o(f),\ f\ne\bar e\}},
 \quad
 (e,f\in \arcs).
\end{equation}
A directed graph is called \emph{strongly connected} if there is a directed path from any vertex to any other. The directed graph associated with $W$ has vertex set $\arcs$ and a directed edge $e\to f$ whenever $W_{ef}=1$; it is strongly connected exactly when $W$ is irreducible.

Write $\uu=(u_e)_{e\in E}$ for the weights on unoriented edges and give each orientation the weight of its underlying edge, so that $u_{\bar e}=u_e$. For a vector $x=(x_i)_{i\in\mathcal I}$ indexed by a finite set $\mathcal I$, we write 
\begin{equation}
D_x:=\diag(x_i)_{i\in\mathcal I}
\end{equation}
for the diagonal matrix with diagonal entries $x_i$. In particular, $D_{\uu}$ is indexed by $\arcs$, with equal entries for the two orientations of each edge.
We then set
\begin{equation}\label{eq:Pintro}
 B(\uu) := WD_{\uu},\qquad
 P_G(\uu) := \det(I-B(\uu)),\qquad
 \zeta_G(\uu) := P_G(\uu)^{-1}.
\end{equation}
The Euler product over primitive reduced cycles converges absolutely to $\zeta_G(\uu)$ whenever $\rho(B(|\uu|))<1$, where $|\uu|$ is taken componentwise and $\rho(A)$ is the largest modulus of an eigenvalue of $A$. The determinant defines $P_G$ as a polynomial and $\zeta_G$ as a rational function for arbitrary weights.

For $0<u_e<1$, define
\begin{equation}\label{eq:hjintro}
 j_{vw}=\frac{u_{vw}}{1-u_{vw}^2}\quad(vw\in E),\qquad
 h_v=1+\sum_{w\sim v}\frac{u_{vw}^2}{1-u_{vw}^2}\quad(v\in V), 
\end{equation}
where $w\sim v$ means that $v$ and $w$ are joined by an edge. Define the symmetric vertex matrix by
\begin{align}
  H_G(\uu)_{vw} := \begin{cases}
    h_v &\quad (v=w) \\
    -j_{vw} &\quad (v\sim w) \\
    0 &\quad (\text{otherwise})
  \end{cases}.
  \label{eq:HG}
\end{align}
%let $H_G(\uu)$ have diagonal entries $h_v$, off-diagonal entries $-j_{vw}$ on edges, and zero otherwise. 
The weighted Ihara--Bass formula gives
\begin{equation}\label{eq:bassintro}
 P_G(\uu)=\prod_{e\in E}(1-u_e^2)\det H_G(\uu).
\end{equation}
The arbitrary-weight form used here is standard; see \citet[Corollary~9]{WatanabeFukumizu2011}. We include a proof to fix conventions.

Our problem is to construct pairs $(\uu',S)$ satisfying
\begin{equation}\label{eq:congruenceintro}
 S H_G(\uu')S=H_G(\uu),\qquad
 S=\diag(a_v)_{v\in V},\quad a_v>0,\quad 0<u'_e<1.
\end{equation}
This is a diagonal \emph{congruence}, not a similarity. It preserves inertia, not individual eigenvalues. Taking determinants gives
\begin{equation}\label{eq:prefactorintro}
 P_G(\uu')=C_G(\uu;\uu',S)^{-1}P_G(\uu),\qquad
 C_G=\prod_v a_v^2\prod_{e\in E}\frac{1-u_e^2}{1-u_e'^2}.
\end{equation}
Thus each solution of \eqref{eq:congruenceintro} gives an identity between the corresponding zeta functions wherever they are finite. The substantive problem is to construct the nontrivial solutions and to determine the range in which they can be classified.

\subsection{Main results}
We now summarize the main results; the objects appearing below are defined in full in \cref{sec:prelim,sec:correspondence,sec:positivebranch}. The first result is a bijection between positive solutions of \eqref{eq:congruenceintro} and solutions of the equations
\begin{equation}\label{eq:cavityintro}
 c_{v\to w}=h_v-\sum_{x\sim v,\ x\ne w}\frac{j_{xv}^2}{c_{x\to v}}.
\end{equation}
These are the scalar Gaussian \emph{cavity} equations. In this correspondence, all $c_{v\to w}$ and the vertex coefficients $s_v$, defined by
\begin{equation}\label{eq:sintro}
 s_v:=h_v-\sum_{x\sim v}\frac{j_{xv}^2}{c_{x\to v}}
 \qquad(v\in V)
\end{equation}
are required to be positive. 
Note that the sum in \eqref{eq:sintro} runs over all neighbours of $v$, whereas the sum on the right-hand side of \eqref{eq:cavityintro} excludes $w$. 
%Thus positivity of all $c_{v\to w}$ does not by itself imply $s_v>0$.

On a tree, removing the interaction across the edge $vw$ separates the branch containing $v$ from the branch containing $w$. Keeping the original diagonal coefficients, $c_{v\to w}$ is the precision at $v$ after integrating out the other variables in its branch. This is the usual cavity construction. On a graph containing cycles, we take $c_{v\to w}$ to be an algebraic unknown in \eqref{eq:cavityintro}. The same equations govern the quadratic coefficients at a fixed point of Gaussian belief propagation (BP) \cite{WeissFreeman2001,MalioutovJohnsonWillsky2006}. \Cref{sec:correspondence} derives the tree recursion and proves the algebraic correspondence on general graphs.

A distinguished reference solution is $c^{(0)}_e=(1-u_e^2)^{-1}$. We obtain a complete classification in the order region $\cav\ge\cav^{(0)}$. Let $\dd_e=1-u_e^2$ and put
\begin{equation}\label{eq:Mintro}
 M_{\uu}=D_{\dd}W^{\mathsf T}D_{\uu^2/\dd},\qquad
 z_e=(1-u_e^2)c_e-1.
\end{equation}
Then \eqref{eq:cavityintro} is equivalent to
\begin{equation}\label{eq:Fintro}
 \zz=M_{\uu}\frac{\zz}{\one+\zz},
\end{equation}
where $\one=(1,\ldots,1)^{\mathsf T}$ is the all-ones vector of the appropriate dimension and division is componentwise.
Assuming that $W$ is irreducible, equation~\eqref{eq:Fintro} has only the zero solution among nonnegative vectors if $\rho(B(\uu^2))\le1$. If $\rho(B(\uu^2))>1$, it has exactly two nonnegative solutions: the zero solution and a unique strictly positive solution.

On the open set
\begin{equation}\label{eq:omega}
 \Omega_+=\{\uu\in(0,1)^E:\rho(B(\uu^2))>1\},
\end{equation}
the positive solution determines the following identity for the weighted Ihara zeta function.

\begin{theorem}[Ihara zeta identity under a change of edge weights]\label{thm:main}
Let $G$ be a finite connected simple graph whose edge adjacency matrix $W$ is irreducible, and let $\uu\in\Omega_+$. Then the system
\begin{equation}\label{eq:mainsystem}
 z_{v\to w}=(1-u_{vw}^2)\sum_{x\sim v,\ x\ne w}
 \frac{u_{xv}^2}{1-u_{xv}^2}\,\frac{z_{x\to v}}{1+z_{x\to v}}
 \qquad(v\to w\in\arcs)
\end{equation}
has exactly one entrywise positive solution $\zz_*(\uu)$. For the corresponding cavity solution $c_{v\to w}=(1+z_{v\to w})/(1-u_{vw}^2)$, let $s_v$ be the vertex coefficients defined in \eqref{eq:sintro}. The transformed weights, these vertex coefficients, and a prefactor are given by
\begin{equation}\label{eq:maindata}
 T_G(\uu)_{vw}=\frac{u_{vw}}{\sqrt{(1+z_{v\to w})(1+z_{w\to v})}},\qquad
 s_v=1+\sum_{x\sim v}\frac{u_{xv}^2}{1-u_{xv}^2}\,\frac{z_{x\to v}}{1+z_{x\to v}},
\end{equation}
\begin{equation}\label{eq:mainC}
 C_G(\uu)=\prod_{v\in V}s_v\prod_{e\in E}\frac{1-u_e^2}{1-T_G(\uu)_e^2},
\end{equation}
all evaluated at $\zz=\zz_*(\uu)$. Then $0<T_G(\uu)_e<u_e<1$, $s_v>1$, $C_G(\uu)>0$, and
\begin{equation}\label{eq:mainintro}
 \zeta_G\bigl(T_G(\uu)\bigr)=C_G(\uu)\,\zeta_G(\uu)
\end{equation}
wherever both sides are finite; equivalently $P_G(T_G(\uu))=C_G(\uu)^{-1}P_G(\uu)$, an identity that also holds at the poles of the zeta functions. The maps $\zz_*$, $T_G$, and $C_G$ are real analytic on $\Omega_+$, and $\rho(B(T_G(\uu)^2))<1$.
\end{theorem}

\Cref{thm:main} is proved as \cref{thm:selected} together with \cref{prop:image}. The identity comes from a positive diagonal congruence $SH_G(T_G(\uu))S=H_G(\uu)$ with $S=\diag\sqrt{s_v}$, and $\zz_*$ is computed by a monotone iteration; no evaluation of either zeta function enters the construction. For uniform weights on regular and biregular graphs, \eqref{eq:mainintro} reduces to the classical reciprocal relations \eqref{eq:regularrelation} and \eqref{eq:biregularp}--\eqref{eq:biregularscales}.
\Cref{ex:diamond} explicitly computes the weight transformation and the resulting functional equation for the diamond graph, which is neither regular nor biregular.

By \cref{thm:main}, the selected transformation sends weights with $\rho(B(\uu^2))>1$ to weights with $\rho(B(T_G(\uu)^2))<1$. It is not a self-map of $\Omega_+$, and we do not claim that this selected branch is an involution. Its inverse congruence is present in the full cavity correspondence, but belongs to a different order region. \Cref{rem:irreducible} characterizes when $W$ is irreducible and shows that $\Omega_+$ is then nonempty.

We also prove that no nontrivial positive diagonal congruence exists when $H_G(\uu)$ is positive definite. Conversely, beyond the nonnegative order region the number of positive cavity solutions is not bounded by two: a chain of $k$ copies of $K_4$, joined by bridges of sufficiently small positive weight, has at least $2^k$ such branches. This is a lower bound for labelled solutions, not a classification of all real or complex branches.

\subsection{Relation to previous work}
\paragraph{Weighted graph zeta functions.}
The determinant formulas developed in \cite{StarkTerras1996,MizunoSato2004,HortonStarkTerras2008,IshikawaMorita2026} for multivariable and weighted zeta functions form the basis of our construction.
Weighted graph zeta functions have also been studied and extended in connection with backtracking parameters \cite{Sato2006,Sato2007}, weighted digraphs and their coverings and products \cite{ChoeEtAl2007}, non-backtracking random walks \cite{Kempton2016}, and discrete-time quantum walks \cite{KonnoSato2012,Ishikawa2022QuantumWalks}.
%We work with scalar edge weights that are symmetric under edge reversal and multiplicative along paths, so the weighted Ihara--Bass formula expresses the zeta function through the vertex matrix $H_G(\uu)$.
We do not introduce a Bartholdi-type backtracking parameter \cite{Bartholdi1999} in this paper.  
%Our question is which pairs of weights $\uu,\uu'$ satisfy $SH_G(\uu')S=H_G(\uu)$ for a positive diagonal matrix $S$.
%Such a congruence yields an explicit prefactor relating the two zeta functions.
%The contribution here is to parameterize these congruences by positive cavity solutions with positive associated vertex coefficients and to construct a specified analytic branch throughout an explicit weight domain.

%Bartholdi-type zeta functions include an additional backtracking parameter \cite{Bartholdi1999,Sato2006}, whose role in functional equations and pole structure is studied by Matsuura and Ohta \cite{MatsuuraOhtaBartholdi2024}.
%Here backtracking is excluded throughout: the transformations act on the independent edge weights of the Ihara zeta function on the same graph.
%For regular and biregular graphs, this construction recovers the familiar reciprocal relations; the diamond graph in \cref{ex:diamond} illustrates how a uniform input can yield nonuniform output weights on an irregular graph.
%
%The matrix weights in the graph-based gauge theories of Matsuura and Ohta \cite{MatsuuraOhta2022,MatsuuraOhtaWilson2022} encode gauge-group representations and lead to Artin--Ihara $L$-functions \cite{MatsuuraOhtaArtinIhara2026}.
%This physical setting motivates the Gaussian-kernel interpretation in \cref{sec:kernels}; the classification and branch construction in this paper concern scalar weights and positive diagonal congruences of their finite vertex matrices.

\paragraph{Matrix weights and graph-based gauge theories.}

Kazakov and Migdal introduced a lattice gauge model with a scalar field in the adjoint representation \cite{KazakovMigdal1993}. Matsuura and Ohta related graph versions of this model to Ihara and Bartholdi zeta functions with matrix-valued edge weights \cite{MatsuuraOhta2022,MatsuuraOhtaWilson2022}. They introduced a fundamental Kazakov--Migdal (FKM) model with fundamental scalar fields \cite{MatsuuraOhtaFKM2023}, studied its phases and duality \cite{MatsuuraOhta2024}, and analyzed functional equations and poles of the Bartholdi zeta function \cite{MatsuuraOhtaBartholdi2024}.

Artin-type $L$-functions for finite graphs and their coverings were developed by Hashimoto \cite{Hashimoto1992} and by Stark and Terras \cite{StarkTerras2000}. In the graph gauge theory, Matsuura and Ohta identify the matrix-weighted zeta functions with Artin--Ihara $L$-functions of graph bundles associated with representations of the gauge group. On a cycle graph, they obtain a random-partition description for the fundamental representation in the large-$N_c$ limit, where $N_c$ is the number of colours \cite{MatsuuraOhtaArtinIhara2026}. The edge weights in our congruence problem are scalar.

\paragraph{Gaussian BP and determinant identities.}
On trees, equations \eqref{eq:cavityintro} follow from exact Gaussian elimination.
On general graphs, the same equations are the fixed-point equations for the quadratic coefficients in Gaussian BP \cite{WeissFreeman2001,MalioutovJohnsonWillsky2006}. In the variational formulation for discrete factor graphs, BP fixed points correspond to stationary points of the Bethe free energy subject to local consistency constraints \cite{YedidiaFreemanWeiss2005}.
Johnson, Chernyak, and Chertkov express a Gaussian determinant as a BP factor times a determinant built from a weighted edge adjacency matrix \cite[Theorem~5]{JohnsonChernyakChertkov2009}.
Their orbit-product arguments impose walk-summability. Chernyak and Chertkov also give a finite loop expansion for determinants without a positive-definiteness requirement, allowing expansion about unstable Gaussian BP solutions \cite{ChernyakChertkov2008}.
Accordingly, merely using an unstable cavity solution or comparing two known determinant factorizations is not, by itself, the issue addressed here.
Our question is which solutions give positive diagonal congruences within the weighted Ihara family, and how a specified branch can be constructed throughout an explicit weight domain.

Watanabe and Fukumizu related the graph zeta function to the Hessian of the Bethe free energy and to the stability of BP fixed points \cite{WatanabeFukumizu2009}. Their subsequent treatment gives the weighted determinant formula and the spectral criterion for positive definiteness used below \cite[Corollary~9 and Theorem~10]{WatanabeFukumizu2011}. The usual scalar specialization of $H_G$ is a rescaled deformed Laplacian, also called a Bethe Hessian in spectral graph methods \cite{SaadeKrzakalaZdeborova2014}. These connections explain why squared edge weights occur in the local stability calculation below.

\paragraph{Nonnegative fixed points and universal covers.}
The substitution $y_i=z_i/(1+z_i)$ converts \eqref{eq:Fintro} into
\[
 0=-y+(I-\diag y)M_{\uu}y.
\]
This is the equilibrium equation of a deterministic susceptible--infected--susceptible (SIS) model with unit recovery rates. The threshold for a unique positive equilibrium is classical \cite{LajmanovichYorke1976}. Van Mieghem, Omic, and Kooij describe the network spectral threshold in the mean-field model \cite{VanMieghemOmicKooij2009}. Further results on deterministic network dynamics appear in \cite{MeiEtAl2017}, and adaptive control of epidemic elimination is studied in \cite{WalshEtAl2025}. Here the coordinates index directed edges. We give a short proof in these variables.

A further relevant line of work concerns the ratio systems of Green functions on universal covers. Angel, Friedman, and Hoory characterize such systems by finite cavity-type equations and a squared-ratio spectral condition \cite{AngelFriedmanHoory2015}. This is closely related to our equations, although their main object is the spectrum of an operator on the universal covering tree rather than a finite-graph weight-congruence relation. We make the algebraic correspondence explicit in \cref{sec:ratios}.

\paragraph{Gaussian-kernel motivation.}
The diagonal-congruence problem also results from rescaling normalized Mehler kernels on the edges and the Gaussian measures at the vertices,
which was inspired by Matsuura and Ohta's Gaussian scalar-field formulations of Kazakov--Migdal models on graphs and their extensions \cite{MatsuuraOhta2022,MatsuuraOhtaFKM2023}.
We retain the scalar, fixed-graph kernel interpretation in \cref{sec:kernels}, but all main identities are proved by finite-dimensional algebra. No assertion about a new physical duality is needed.

Taken together, these works supply the determinant and fixed-point tools. The results below concern their specific application to positive diagonal congruences of weighted Ihara matrices, including the order-restricted classification, the selected analytic transformation, rigidity, and the multiplicity construction. We do not claim priority for the separate determinant or nonnegative fixed-point theorems, and do not infer global uniqueness from uniqueness in the nonnegative order region.

\subsection{Outline of the paper}
\Cref{sec:prelim} fixes conventions, records the Euler product, and proves the weighted Ihara--Bass formula (\cref{prop:bass}).
\Cref{sec:correspondence} derives the cavity equations from exact Gaussian elimination on trees, proves the bijection between positive cavity solutions with positive vertex coefficients $s_v$ and positive diagonal congruences (\cref{thm:correspondence}), and rewrites the system relative to the reference solution as the fixed-point problem \eqref{eq:weightedfixed}.
\Cref{sec:positivebranch} proves the main theorem: \cref{lem:fixed,lem:iteration} establish the nonnegative fixed-point classification and monotone iteration, \cref{thm:selected} classifies cavity solutions in the nonnegative order region and constructs the real-analytic weight transformation on $\Omega_+$, \cref{prop:image} shows that its image is subcritical for squared weights, and \cref{sec:ratios} delimits the relation to universal-cover ratio systems.
\Cref{sec:convergence} identifies the positive-definite region with the unsquared spectral condition (\cref{thm:positivity}), discusses convergent Euler products and Gaussian integrals, and proves the rigidity theorem (\cref{thm:rigidity}).
\Cref{sec:examples} classifies uniform-to-uniform congruences and recovers the regular and biregular reciprocal relations.
The diamond graph in \cref{sec:diamondexample} illustrates a uniform input with a nonuniform output, while \cref{sec:vertexinputs} treats vertex-factorized weights and ends with the $K_4$ example in \cref{sec:inhomogeneous}, whose output is not vertex-factorized.
\Cref{sec:manybranches} constructs at least $2^k$ admissible cavity branches by weakly linking $k$ copies of $K_4$ (\cref{thm:manybranches}), demonstrating multiplicity beyond the nonnegative order region; \cref{sec:fourbranches} gives a numerical example with two linked copies.
\Cref{sec:discussion} is devoted to the conclusion and discussion.
\Cref{sec:kernels} explains the Gaussian-kernel motivation and gives a signed-weight example showing why positivity of the edge weights is essential for rigidity.
\Cref{sec:computations} gives pseudocode, input data, and exact checks for reproducing the examples.

\section{Weighted Ihara determinants and conventions}\label{sec:prelim}
All inequalities between vectors or matrices, except positive definiteness, are entrywise. Products over $E$ use each unoriented edge once; matrices on $\arcs$ use both orientations. In particular, $D_{\uu^2}$ has $2m$ entries, whereas $\prod_{e\in E}(1-u_e^2)$ has $m$ factors. For a weight vector $\uu$, the symbol $\uu^2$ always denotes the entrywise square $(u_e^2)_{e\in\arcs}$; thus $B(\uu^2)=WD_{\uu^2}$, which must not be confused with the matrix square $B(\uu)^2$. We write $x>0$ for strict positivity of every component and $H\succ0$ for a positive-definite real symmetric matrix.

\subsection{Euler products and their determinant representation}
A cycle $C=(e_1,\ldots,e_\ell)$ is reduced exactly when $W_{e_i e_{i+1}}=1$ for all $i=1,\cdots \ell$ modulo $\ell$, so immediate reversals are forbidden also at the cyclic join.
Two cycles are equivalent if they differ by a cyclic shift. For $C=(e_1,\ldots,e_\ell)$, denote its equivalence class by
\[
 [C]=\{(e_k,\ldots,e_\ell,e_1,\ldots,e_{k-1}):1\le k\le\ell\}.
\]
The Euler product identifies these cyclic shifts, but does not identify $C$ with the reversed cycle $(\bar e_\ell,\ldots,\bar e_1)$. Write $u(C)=\prod_{i=1}^\ell u_{e_i}$, with each traversal counted separately; this weight is invariant under cyclic shifts. On the domain $\rho(B(|\uu|))<1$, the following Euler product converges absolutely:
\begin{equation}\label{eq:euler}
 \zeta_G(\uu) = \prod_{[C]\text{ primitive}}(1-u(C))^{-1}.
\end{equation}
Here the product runs over equivalence classes of primitive reduced cycles. Indeed,
\begin{equation}\label{eq:traces}
 -\log\det(I-B(\uu))=\sum_{k\ge1}\frac{\Tr B(\uu)^k}{k}.
\end{equation}
A primitive reduced cycle of length $\ell$ repeated $r$ times contributes $\ell$ choices of starting arc to the trace at $k=r\ell$. Its total contribution is therefore $u(C)^r/r$, which gives \eqref{eq:euler}. The identity also holds formally in the independent edge variables.

\begin{proposition}[Weighted Ihara--Bass formula]\label{prop:bass}
If $u_{\bar e}=u_e$ and $u_e^2\ne1$, then
\begin{equation}\label{eq:bass}
 P_G(\uu)=\left(\prod_{e\in E}(1-u_e^2)\right)\det H_G(\uu).
\end{equation}
This is an identity of rational expressions whose left-hand side is polynomial.
\end{proposition}
\begin{proof}
Let $O_{ve}=\mathbf1_{\{v=o(e)\}}$ and $T_{ve}=\mathbf1_{\{v=t(e)\}}$, and let $J_{ef}=\mathbf1_{\{f=\bar e\}}$. Then $W=T^{\mathsf T}O-J$. 
Since edge symmetry gives $JD_{\uu}=D_{\uu}J$ and $J^2=I$, we have $(I+JD_{\uu})(I-JD_{\uu})=I-D_{\uu}^2$. The assumption $u_e^2\ne1$ makes $I-D_{\uu}^2$ invertible, so
\[
 (I+JD_{\uu})^{-1}=(I-JD_{\uu})(I-D_{\uu}^2)^{-1}.
\]
Ordering the arcs in pairs $e,\bar e$ makes $I+JD_{\uu}$ block diagonal, with one block $\begin{psmallmatrix}1&u_e\\u_e&1\end{psmallmatrix}$ for each unoriented edge. Taking the product of the block determinants gives
\[
 \det(I+JD_{\uu})=\prod_{e\in E}(1-u_e^2).
\]
The determinant identity $\det(I-XY)=\det(I-YX)$ yields
\begin{align*}
 P_G(\uu)
 &=\det(I+JD_{\uu}-T^{\mathsf T}O D_{\uu})\\
 &=\det(I+JD_{\uu})\det\bigl(I-O D_{\uu}(I+JD_{\uu})^{-1}T^{\mathsf T}\bigr).
\end{align*}
The second determinant is $\det H_G(\uu)$: its off-diagonal term is $-u_{vw}/(1-u_{vw}^2)$, and the term involving $J$ adds $\sum_{w\sim v}u_{vw}^2/(1-u_{vw}^2)$ to its $v$th diagonal. This proves \eqref{eq:bass}.
\end{proof}

\begin{remark}
The formula is standard \cite{WatanabeFukumizu2009,WatanabeFukumizu2011}. Positivity of $H_G(\uu)$ is not required for the determinant identity. This distinction will be essential: the nontrivial transformations selected below do not lie in the positive-definite region.
\end{remark}

\section{Cavity equations and positive diagonal congruence}\label{sec:correspondence}

%For fixed $\uu\in(0,1)^E$, we define $h_v$ and $j_{vw}$ by \eqref{eq:hjintro} and $H_G$ by \eqref{eq:HG}.
\subsection{Exact Gaussian elimination on trees}\label{sec:treeelimination}

The equations studied in this section arise from Gaussian elimination, and we first recall why. When $H_G(\uu)\succ0$, consider the Gaussian density on $\CC^V$ proportional to $e^{-X^\dagger H_G(\uu)X}$, where $X=(X_v)_{v\in V}$.

Suppose first that $G$ is a tree. For an arc $v\to w$, delete the edge $\{v,w\}$ and let $V_{v\mid w}$ be the vertex set of the component containing $v$. Let $H^{(v\mid w)}$ be the corresponding principal submatrix of $H_G(\uu)$, retaining the original diagonal entries $h_a$, and write $Y=(X_a)_{a\in V_{v\mid w}}$. Keep $X_v$ fixed and integrate all other variables in this component. Define $c_{v\to w}$ by the resulting Gaussian factor:
\begin{equation}\label{eq:cavityintegral}
 \int_{\CC^{V_{v\mid w}\setminus\{v\}}}
 e^{-Y^\dagger H^{(v\mid w)}Y}
 \prod_{a\in V_{v\mid w}\setminus\{v\}}\frac{d^2X_a}{\pi}
 =A_{v\to w}e^{-c_{v\to w}|X_v|^2}.
\end{equation}
Here $d^2X_a=d(\Re X_a)d(\Im X_a)$ and $A_{v\to w}>0$ is independent of $X_v$. Thus $c_{v\to w}>0$ is the coefficient of $|X_v|^2$ in the effective quadratic form after the $w$-side branch has been excluded. Equivalently, it is the precision of $X_v$ in this cavity model. If $v$ is a leaf with neighbour $w$, the component contains only $v$, so $A_{v\to w}=1$ and $c_{v\to w}=h_v$.
Note that the coefficients are indexed by arcs: $c_{v\to w}$ refers to the component containing $v$, and $c_{w\to v}$ to the component containing $w$.

To derive the equation for $c_{v\to w}$, keep $X_v$ fixed and consider a neighbour $x\ne w$. Applying \eqref{eq:cavityintegral} to the component $V_{x\mid v}$ first integrates all variables in that branch except $X_x$, leaving the factor $A_{x\to v}e^{-c_{x\to v}|X_x|^2}$. Including the coupling between $x$ and $v$ and completing the square then gives
\begin{equation}\label{eq:cavitytransfer}
 \int_{\CC}
 e^{-c_{x\to v}|X_x|^2+j_{xv}(\overline X_xX_v+\overline X_vX_x)}
 \frac{d^2X_x}{\pi}
 =\frac1{c_{x\to v}}
   \exp\!\left(\frac{j_{xv}^2}{c_{x\to v}}|X_v|^2\right).
\end{equation}
The full $x$-branch therefore contributes the prefactor $A_{x\to v}/c_{x\to v}$ and reduces the coefficient of $|X_v|^2$ by $j_{xv}^2/c_{x\to v}$. For an original leaf $x$, this is the replacement $h_v\mapsto h_v-j_{xv}^2/h_x$. Since $G$ is a tree, the branches indexed by $x\sim v$, $x\ne w$, are disjoint and have no couplings to one another. Integrating all these branches therefore leaves the coefficient
\begin{equation}\label{eq:cavityupdate}
 c_{v\to w}=[\mathcal Q(\cav)]_{v\to w}:=h_v-\sum_{x\sim v,\ x\ne w}\frac{j_{xv}^2}{c_{x\to v}}.
\end{equation}
This defines a rational map $\mathcal Q$ on vectors with nonzero entries, and the recursion can be written as $\cav=\mathcal Q(\cav)$: $\cav$ is given by the fixed point of the mapping ${\mathcal Q}$.
On a tree, successive leaf integrations terminate and compute $\det H_G(\uu)$ exactly.

For the particular coefficients \eqref{eq:hjintro}, the recursion \eqref{eq:cavityupdate} implies that on a tree
\begin{equation}\label{eq:treecavity}
 c_{v\to w}=c_{w\to v}=\frac1{1-u_{vw}^2}.
\end{equation}
The leaf values above have this form. If the incoming coefficients have this form, each term $j_{xv}^2/c_{x\to v}=u_{xv}^2/(1-u_{xv}^2)$ cancels the corresponding contribution to $h_v$, leaving $1+u_{vw}^2/(1-u_{vw}^2)=(1-u_{vw}^2)^{-1}$. This proves \eqref{eq:treecavity} by induction from the leaves.

\subsection{The cavity correspondence on general graphs}\label{sec:cavitycorrespondence}
We now consider a general finite simple graph $G$. On a tree, deleting an edge separates the graph into two components, allowing the independent branch integrations in \eqref{eq:cavityintegral}. 
If an edge lies on a cycle, however, deleting it leaves its endpoints connected by another path, so this decomposition into independent branches is no longer available.

The Gaussian integrations themselves remain well-defined when $H_G(\uu)\succ0$, but elimination also changes the couplings between the remaining variables. First consider the triangle $v-x-y-v$. Unlike a leaf, the vertex $x$ has two neighbours, $v$ and $y$, that remain after $X_x$ is integrated out. Keeping $X_v$ and $X_y$ fixed, the factor involving $X_x$ gives
\begin{align*}
 \int_{\CC}
 \exp\!\left[-h_x|X_x|^2
 +2\Re\!\left(\overline X_x(j_{xv}X_v+j_{xy}X_y)\right)\right]
 \frac{d^2X_x}{\pi}
 =\frac1{h_x}
 \exp\!\left(\frac{|j_{xv}X_v+j_{xy}X_y|^2}{h_x}\right).
\end{align*}
Combining this factor with the original quadratic form on $v$ and $y$ gives the effective coefficients
\begin{align*}
 \widetilde h_v=h_v-\frac{j_{xv}^2}{h_x},\qquad 
 \widetilde h_y=h_y-\frac{j_{xy}^2}{h_x},\qquad
 \widetilde j_{vy}=j_{vy}+\frac{j_{xv}j_{xy}}{h_x}.
\end{align*}
The edge $\{v,y\}$ is already present in the triangle, but its coupling must be updated along with the diagonal coefficients. A subsequent integration of $X_y$ uses $\widetilde h_y$ and $\widetilde j_{vy}$, so the effects of eliminating $x$ and $y$ cannot be computed as independent branch contributions with the original couplings held fixed.

For a four-cycle $v-x-y-w-v$ with no diagonal edges, elimination produces a phenomenon that does not occur in the triangle.
Integrating $X_x$ while keeping $X_v$, $X_y$, and $X_w$ fixed produces the same correction to the quadratic form, including the cross term
\[
 -\frac{j_{xv}j_{xy}}{h_x}
 \bigl(\overline X_vX_y+\overline X_yX_v\bigr).
\]
Here $\{v,y\}$ is not an edge of the original graph, so this term creates a new edge in the graph of effective couplings. 
On a graph containing cycles, exact elimination can modify existing couplings or create new edges.
The triangle illustrates the first effect and the four-cycle the second.
Successive leaf integrations on a tree always eliminate a vertex with only one remaining neighbour.
Each integration therefore changes only the diagonal coefficient of that neighbour.
These distinctions explain why exact elimination on a general graph is not described by the scalar recursion \eqref{eq:cavityupdate} with the original graph and couplings held fixed.

From this point on, we call a vector $\cav$ satisfying \eqref{eq:cavityupdate} for every arc a \emph{cavity solution}. Its entries are algebraic unknowns and, on graphs containing simple cycles, need not coincide with coefficients obtained by exact partial Gaussian integration. 
In addition, unlike the tree values in \eqref{eq:treecavity}, these entries need not satisfy $c_{v\to w}=c_{w\to v}$. 
This algebraic system does not assume positive definiteness, convergence of Gaussian integrals, or convergence of message iterations. For context, \eqref{eq:cavityupdate} is also the fixed-point equation for the quadratic-coefficient update in Gaussian BP \cite{WeissFreeman2001,MalioutovJohnsonWillsky2006}.
\begin{definition}[Admissible cavity solution]\label[definition]{def:admissible}
For a cavity solution $\cav$, write the vertex coefficients \eqref{eq:sintro} again as
\begin{equation}\label{eq:sdef}
 s_v(\cav):=h_v-\sum_{x\sim v}\frac{j_{xv}^2}{c_{x\to v}}
 \quad(v\in V).
\end{equation}
We call $\cav$ an \emph{admissible cavity solution} if $c_{v\to w}>0$ for every arc $v\to w$ and $s_v(\cav)>0$ for every vertex $v$. The condition $\cav>0$ does not imply $s_v(\cav)>0$ for every vertex $v$; admissibility requires both conditions.
\end{definition}
For positive cavity coefficients, the vertex condition has a direct meaning. The reconstruction in \eqref{eq:fromc} gives $u'_{vw}>0$, and \eqref{eq:csidentity} gives $1-(u'_{vw})^2=s_v/c_{v\to w}$. Thus $s_v>0$ at every vertex is exactly the condition that every reconstructed edge weight is less than one.

\begin{theorem}[Cavity parameterization]\label{thm:correspondence}
There is a bijection between the admissible cavity solutions of \cref{def:admissible} and pairs $(\uu',S)$ with $\uu'\in(0,1)^E$ and $S=\diag(a_v)>0$ satisfying $SH_G(\uu')S=H_G(\uu)$. In the forward direction,
\begin{equation}\label{eq:fromc}
 a_v^2=s_v(\cav),\qquad
 u'_{vw}=\frac{j_{vw}}{\sqrt{c_{v\to w}c_{w\to v}}}.
\end{equation}
In the reverse direction,
\begin{equation}\label{eq:toc}
 c_{v\to w}=\frac{a_v^2}{1-u_{vw}'^{\,2}}.
\end{equation}
For every such pair, the determinant relation \eqref{eq:prefactorintro} holds, including when either determinant vanishes.
\end{theorem}
\begin{proof}
Subtracting the full vertex sum in \eqref{eq:sdef} from the sum with one neighbor omitted gives, for each edge $vw$,
\begin{equation}\label{eq:csidentity}
 s_v=c_{v\to w}-\frac{j_{vw}^2}{c_{w\to v}}
 =c_{v\to w}(1-u_{vw}'^{\,2}).
\end{equation}
Because $s_v$ and $c_{v\to w}$ are both positive, the value in \eqref{eq:fromc} satisfies $0<u'_{vw}<1$. 
Applying \eqref{eq:csidentity} at both ends $v$ and $w$ and taking positive square roots yields
\begin{equation}\label{eq:crossmatch}
 a_va_w\frac{u'_{vw}}{1-u_{vw}'^{\,2}}=j_{vw}.
\end{equation}
Moreover,
\[
 s_v\frac{u_{vw}'^{\,2}}{1-u_{vw}'^{\,2}}
 =\frac{j_{vw}^2}{c_{w\to v}}.
\]
Summing over $w$ and using \eqref{eq:sdef} gives
\begin{equation}\label{eq:diagmatch}
 a_v^2\left(1+\sum_{w\sim v}\frac{u_{vw}'^{\,2}}{1-u_{vw}'^{\,2}}\right)=h_v.
\end{equation}
Equations \eqref{eq:crossmatch} and \eqref{eq:diagmatch} are exactly the off-diagonal and diagonal congruence conditions.

Conversely, assume the congruence and define $c$ by \eqref{eq:toc}. The off-diagonal identity \eqref{eq:crossmatch} then gives
\[
 \frac{j_{xv}^2}{c_{x\to v}}
 =a_v^2\frac{u_{xv}'^{\,2}}{1-u_{xv}'^{\,2}}.
\]
Subtracting these expressions from \eqref{eq:diagmatch} gives $s_v=a_v^2>0$. Omitting the $w$th term instead gives
\[
 h_v-\sum_{x\sim v,\ x\ne w}\frac{j_{xv}^2}{c_{x\to v}}
 =a_v^2\left(1+\frac{u_{vw}'^{\,2}}{1-u_{vw}'^{\,2}}\right)
 =c_{v\to w}.
\]
The constructions are inverse by \eqref{eq:csidentity}. Finally, taking determinants in the congruence and applying \cref{prop:bass} proves the stated polynomial relation, without division by $P_G$.
\end{proof}

\subsection{The explicit reference solution and the order region}
The reference solution
\begin{equation}\label{eq:baseline}
 c_e^{(0)}=\frac1{1-u_e^2}
\end{equation}
has $s_v=1$ and gives the identity transformation. This follows since $j_e^2/c_{\bar e}^{(0)}=u_e^2/(1-u_e^2)$.

Introduce $z_e=(1-u_e^2)c_e-1$. Subtracting the reference cavity equation from \eqref{eq:cavityupdate} gives
\begin{equation}\label{eq:weightedfixed}
 z_{v\to w}=(1-u_{vw}^2)
 \sum_{x\sim v,\ x\ne w}\frac{u_{xv}^2}{1-u_{xv}^2}
 \frac{z_{x\to v}}{1+z_{x\to v}}.
\end{equation}
With $M_{\uu}$ defined in \eqref{eq:Mintro}, we set
\begin{equation}\label{eq:Fz}
 F_{\uu}(\zz)=M_{\uu}\frac{\zz}{\one+\zz}
 \qquad(z_e>-1\text{ for all }e\in\arcs),
\end{equation}
where division is componentwise. Then \eqref{eq:weightedfixed} is the fixed-point equation $\zz=F_{\uu}(\zz)$ and is therefore equivalent to \eqref{eq:Fintro}. For every nonnegative solution, the positivity condition on the vertices is automatic because
\begin{equation}\label{eq:sfromz}
 s_v=1+\sum_{x\sim v}\frac{u_{xv}^2}{1-u_{xv}^2}
                    \frac{z_{x\to v}}{1+z_{x\to v}}\ge1.
\end{equation}
The transformed weights are
\begin{equation}\label{eq:ufromz}
 u'_{vw}=\frac{u_{vw}}{\sqrt{(1+z_{v\to w})(1+z_{w\to v})}}.
\end{equation}
In particular, the order condition $\cav\ge\cav^{(0)}$ is equivalent to $\zz\ge0$. It is a condition on directed cavity coefficients, not merely the condition $\uu'\le\uu$ on undirected weights.

\section{Construction and uniqueness of the nonnegative branch}\label{sec:positivebranch}
\subsection{A standard nonnegative fixed-point theorem}
The fixed-point classification below is the endemic-equilibrium threshold theorem for a deterministic network SIS system in different coordinates \cite{LajmanovichYorke1976,VanMieghemOmicKooij2009,MeiEtAl2017}. We give separate proofs for the classification and for the monotone iteration and stability properties used here.

\begin{lemma}[Nonnegative fixed points]\label[lemma]{lem:fixed}
Let $M$ be a nonnegative irreducible $N\times N$ matrix with $N\ge2$. For $z=(z_1,\ldots,z_N)^{\mathsf T}\in\RR_{\ge0}^N$, define $F$ componentwise by
\[
 F_i(z):=\sum_{j=1}^{N}M_{ij}\frac{z_j}{1+z_j}
 \qquad(i=1,\ldots,N).
\]
If $\rho(M)\le1$, the only fixed point of $F$ is $z=0$. If $\rho(M)>1$, $F$ has a unique positive fixed point $z_*$.
\end{lemma}
\begin{proof}
First, let $z\ge0$ be a fixed point. If $z_i=0$ for some $i$, then
\[
 0=F_i(z)=\sum_{j=1}^{N}M_{ij}\frac{z_j}{1+z_j}.
\]
Every summand is nonnegative, so $z_j=0$ whenever $M_{ij}>0$. Irreducibility propagates this zero condition to every component, giving $z=0$. Hence a nonzero nonnegative fixed point has no zero components and is therefore strictly positive.

Let $\ell>0$ be a left Perron vector of $M$. The corresponding eigenvalue is $\rho(M)$ by definition. If $z>0$ is fixed, then
\[
 \ell^{\mathsf T}z=\rho(M)\ell^{\mathsf T}\frac{z}{\one+z}
 <\rho(M)\ell^{\mathsf T}z,
\]
so $\rho(M)>1$ is necessary.
This also means that $z=0$ is the unique fixed point for $\rho(M)\le 1$. 

For sufficiency, let $\lambda:=\rho(M)>1$, choose $r>0$ with $Mr=\lambda r$ and $\max_i r_i=1$, and define $\alpha=(\lambda-1)/2$. Since $r_i\le1$,
\[
 F(\alpha r)\ge\frac{\alpha}{1+\alpha}Mr
 =\frac{\lambda}{1+\alpha}\alpha r\ge\alpha r.
\]
Also $F(z)\le M\one$ for every $z\ge0$, and
$\alpha r\le\lambda r=Mr\le M\one$.
Since $F$ preserves the componentwise order, the iteration starting at $\alpha r$ is nondecreasing and bounded above by $M\one$. Its limit is a positive fixed point by continuity of $F$.

For uniqueness, let $z,w>0$ be fixed points. If they are different, interchange them if necessary to arrange $\theta:=\max_i z_i/w_i>1$. Then $z\le\theta w$, with equality in at least one component since $N$ is finite. Every row of $M$ is nonzero, because $M$ is irreducible with $N\ge2$, and
\[
 \frac{\theta w_i}{1+\theta w_i}<\theta\frac{w_i}{1+w_i},
\]
so $F(\theta w)<\theta F(w)$. At a component where $z_i=\theta w_i$, monotonicity gives $z_i=F_i(z)\le F_i(\theta w)<\theta w_i$, which is a contradiction. This proves uniqueness.
\end{proof}

\begin{lemma}[Monotone iteration and stability]\label[lemma]{lem:iteration}
Let $M$ and $F$ be as in \cref{lem:fixed}, assume $\lambda:=\rho(M)>1$, and let $z_*>0$ be the unique nonzero nonnegative fixed point. Choose $r>0$ with $Mr=\lambda r$ and $\max_i r_i=1$. Then the sequences $(L_k)$ and $(U_k)$ ($k=0,1,2,\cdots$) defined by 
\begin{equation}\label{eq:brackets}
 L_{k+1}=F(L_k),\quad U_{k+1}=F(U_k),
\end{equation}
with 
\begin{equation}
 L_0=\frac{\lambda-1}{2}r,\qquad U_0=M\one
\end{equation}
satisfy
\[
 L_0\le L_1\le L_2\le\cdots\le z_*\le\cdots\le U_2\le U_1\le U_0
\]
componentwise and both sequences converge to $z_*$. 
The iteration thus brackets the fixed point between a nondecreasing lower sequence and a nonincreasing upper sequence. 
Moreover, 
with vector powers interpreted componentwise, 
we have
\begin{equation}\label{eq:stablefixed}
 \rho\bigl(DF(z_*)\bigr)=\rho\bigl(MD_{(\one+z_*)^{-2}}\bigr)<1, 
\end{equation} 
where $DF(z):=(\partial F_i(z)/\partial z_j)_{i,j=1}^N$ is the Jacobian matrix of $F$.
\end{lemma}
\begin{proof}
The bounds in the proof of \cref{lem:fixed} give $L_0\le F(L_0)$, $L_0\le U_0$, and $F(z)\le U_0$ for every $z\ge0$. Since $F$ preserves the componentwise order, $(L_k)$ is nondecreasing, $(U_k)$ is nonincreasing, and $L_k\le U_k$ for every $k$. Their limits are fixed points by continuity, and the lower limit is at least $L_0>0$. Uniqueness in \cref{lem:fixed} therefore implies that both limits equal $z_*$, proving the stated bounds and convergence.

For stability, differentiation gives $DF(z_*)=MD_{(\one+z_*)^{-2}}$. Since every row of $M$ is nonzero,
\[
 DF(z_*)z_*=M\frac{z_*}{(\one+z_*)^2}
 <M\frac{z_*}{\one+z_*}=z_*.
\]
Similarity by $D_{z_*}$ bounds the spectral radius by the largest row sum, which is strictly below one. This proves \eqref{eq:stablefixed}.
\end{proof}

\begin{remark}[SIS equilibria and discrete iteration]
The componentwise change of variables $y_i=z_i/(1+z_i)$, with inverse $z_i=y_i/(1-y_i)$, is a bijection between $[0,\infty)^N$ and $[0,1)^N$. It transforms the fixed-point condition $z=F(z)$ into $0=-y+(I-\diag y)My$. Thus the fixed points of $F$ correspond exactly to the equilibria of the deterministic network SIS differential equation with unit recovery rates \cite{MeiEtAl2017},
\[
 \dot y=-y+(I-\diag y)My.
\]
Under the same change of variables, the discrete iteration $z^{(k+1)}=F(z^{(k)})$ in \cref{lem:iteration} instead becomes
\[
 y^{(k+1)}=\frac{My^{(k)}}{\one+My^{(k)}},
\]
with componentwise division. The stationary states therefore correspond exactly, whereas the discrete iteration and the continuous-time SIS equation have different evolution rules. This correspondence transfers a deterministic equilibrium theorem; it does not identify the discrete iterates with the continuous-time SIS trajectories. The correspondence applies only to $z\ge0$ and gives no conclusion about fixed points with negative $z_i$.
\end{remark}

\subsection{Application to the weighted graph}
From the definition of $M_{\uu}$ in \eqref{eq:Mintro}, using $D_{\uu^2/\dd}D_{\dd}=D_{\uu^2}$, we obtain
\begin{equation}\label{eq:Mspectrum}
 D_{\dd}^{-1}M_{\uu}D_{\dd}=W^{\mathsf T}D_{\uu^2},
 \qquad \rho(M_{\uu})=\rho(B(\uu^2)).
\end{equation}
Since positive diagonal factors do not change the zero pattern, $M_{\uu}$ is irreducible exactly when $W$ is. 
The spectral-radius equality follows by transposition and the similarity of $D_{\uu^2}W$ to $WD_{\uu^2}$.

\begin{remark}[Irreducibility and nonemptiness of $\Omega_+$]\label{rem:irreducible}
For a finite connected simple graph, $W$ is irreducible exactly when every vertex has degree at least two and the graph is not a cycle graph. Indeed, an arc into a vertex of degree one has no non-backtracking successor and creates a zero row; on a cycle graph the two orientation classes do not communicate; in all remaining cases any directed edge can be joined to any other by a non-backtracking path \cite[Theorem~11.10(3) and Lemma~11.11]{Terras2011}. An irreducible $W$ has all row sums at least one and, not being a permutation matrix, not all equal to one; for an irreducible nonnegative matrix this forces $\rho(W)>1$. Since $\rho(B(\uu^2))\to\rho(W)$ as all $u_e\to1$, the set $\Omega_+$ of \eqref{eq:omega} is nonempty whenever $W$ is irreducible. For a cycle graph with $n$ vertices, by contrast, $W$ splits into its two orientation classes, $\rho(B(\uu^2))=\bigl(\prod_{e\in E}u_e^2\bigr)^{1/n}<1$ on each, and no $\uu\in(0,1)^E$ lies in $\Omega_+$.
\end{remark}

\emph{The next theorem is the central result of the paper}; together with \cref{prop:image} below, it proves \cref{thm:main}.

\begin{theorem}[Selected weight transformation]\label{thm:selected}
Assume $W$ irreducible. In the order region $c_e\ge(1-u_e^2)^{-1}$, the reference solution $\cav^{(0)}$ defined in \eqref{eq:baseline} is the only cavity solution if $\rho(B(\uu^2))\le1$. If $\uu\in\Omega_+$, this order region contains exactly two cavity solutions: $\cav^{(0)}$ and a solution $\cav$ satisfying $c_e>c_e^{(0)}$ for every arc $e\in\arcs$.

Let $z_*(\uu)$ be the coordinate vector of this second solution. Equations \eqref{eq:sfromz} and \eqref{eq:ufromz} define $S_G(\uu)=\diag\sqrt{s_v}$ and $T_G(\uu)=\uu'$ which satisfy 
\begin{gather}\label{eq:selectedcongruence}
 S_G(\uu)H_G(T_G(\uu))S_G(\uu)=H_G(\uu),\\
 0<T_G(\uu)_e<u_e<1,\qquad s_v>1.\label{eq:selectedbounds}
\end{gather}
The maps $z_*$, $S_G$, $T_G$ and
\begin{equation}\label{eq:Cselected}
 C_G(\uu)=\prod_vs_v\prod_{e\in E}\frac{1-u_e^2}{1-T_G(\uu)_e^2}
\end{equation}
are real analytic on $\Omega_+$. They satisfy
\begin{equation}\label{eq:mainP}
 P_G(T_G(\uu))=C_G(\uu)^{-1}P_G(\uu),
\end{equation}
which is equivalent to 
%Where the reciprocals are finite, equivalently
\begin{equation}\label{eq:mainzeta}
 \zeta_G(T_G(\uu))=C_G(\uu)\zeta_G(\uu)
\end{equation}
when the reciprocals are finite. 
\end{theorem}
\begin{proof}

With $z_e=(1-u_e^2)c_e-1$, the cavity equation \eqref{eq:weightedfixed} takes the form $\zz=F_{\uu}(\zz)$ in \eqref{eq:Fz}. This is the map of \cref{lem:fixed} with $M=M_{\uu}$. Since $c_e=(1+z_e)/(1-u_e^2)$ and $c_e^{(0)}=(1-u_e^2)^{-1}$, we have
\[
 \cav\ge\cav^{(0)}\ \Longleftrightarrow\ \zz\ge0,
 \qquad
 \cav=\cav^{(0)}\ \Longleftrightarrow\ \zz=0.
\]
%By \eqref{eq:weightedfixed}, the cavity equations are equivalent to $\zz=F_{\uu}(\zz)$. 
Thus cavity solutions in this order region correspond bijectively to nonnegative fixed points of $F_{\uu}$.

The matrix $M_{\uu}$ is nonnegative and irreducible because $W$ is irreducible, and \eqref{eq:Mspectrum} gives $\rho(M_{\uu})=\rho(B(\uu^2))$. If $\rho(B(\uu^2))\le1$, \cref{lem:fixed} shows that $\zz=0$ is the only nonnegative fixed point, so $\cav^{(0)}$ is the only cavity solution in the order region. If $\uu\in\Omega_+$, then $\rho(M_{\uu})>1$, and \cref{lem:fixed} and its proof show that the nonnegative fixed points are exactly $0$ and a unique positive vector $z_*(\uu)$. These correspond to $\cav^{(0)}$ and a second cavity solution satisfying
\[
 c_e-c_e^{(0)}=\frac{(z_*(\uu))_e}{1-u_e^2}>0
 \qquad(e\in\arcs).
\]
This proves the claimed classification.

For every nonnegative fixed point, \eqref{eq:sfromz} gives $s_v\ge1$, so the corresponding cavity solution is admissible. For the positive fixed point $z_*(\uu)$, every summand in \eqref{eq:sfromz} is positive and every vertex has a neighbour, giving $s_v>1$. Also, the denominator in \eqref{eq:ufromz} is greater than one on every edge, so $0<u'_{vw}<u_{vw}<1$. Applying \cref{thm:correspondence} to this admissible cavity solution gives the stated congruence and determinant identities, and hence the zeta identity wherever the reciprocals are finite.

The map $F_{\uu}(z)$ is real analytic for $0<u_e<1$ and $z_e>-1$. By \eqref{eq:stablefixed}, $I-D_zF_{\uu}(z_*)$ is invertible. The analytic implicit-function theorem gives a local analytic fixed-point branch at each $\uu\in\Omega_+$. Positivity and uniqueness make these local branches agree on overlaps. Positive square roots in \eqref{eq:ufromz} and $S_G$ preserve analyticity.
\end{proof}

The selected branch is constructive through \eqref{eq:brackets}. In exact arithmetic these iterates bracket the solution. A floating-point implementation is an illustration, not a validated enclosure unless rounding is controlled. For numerical stability one may use $\eta_e=z_e/(1-u_e^2)=c_e-c_e^{(0)}$, in which the iteration becomes
\begin{equation}\label{eq:eta}
 \eta^{(k+1)}_{v\to w}=\sum_{x\sim v,\ x\ne w}
 \frac{u_{xv}^2\eta^{(k)}_{x\to v}}{1+(1-u_{xv}^2)\eta^{(k)}_{x\to v}}.
\end{equation}

\begin{remark}[Boundary behaviour]\label{rem:boundary}
Let $\uu_0\in(0,1)^E$ satisfy $\rho(B(\uu_0^2))=1$. Along any sequence in $\Omega_+$ converging to $\uu_0$, the positive solutions are bounded above by $z_*(\uu)\le M_{\uu}\one$, and every subsequential limit is a nonnegative fixed point of $F_{\uu_0}$, hence zero by \cref{lem:fixed}. Therefore $z_*(\uu)\to0$, and $T_G$, $S_G$, and $C_G$ extend continuously to such boundary points by $T_G(\uu_0)=\uu_0$, $S_G(\uu_0)=I$, and $C_G(\uu_0)=1$. The selected branch thus emerges from the reference solution at the threshold.
\end{remark}

\subsection{The image crosses the squared-weight threshold}\label{sec:thresholdimage}
Recall from \eqref{eq:cavityupdate} that the cavity update map $\mathcal Q$ is given componentwise by
\[
 [\mathcal Q(\cav)]_{v\to w}
 =h_v-\sum_{x\sim v,\ x\ne w}\frac{j_{xv}^2}{c_{x\to v}}
 \qquad(v\to w\in\arcs).
\]
\begin{proposition}\label[proposition]{prop:image}
\begin{enumerate}[label=(\arabic*)]
\item 
For any admissible cavity solution $\cav$ in the sense of \cref{def:admissible}, let $\uu'$ be the transformed weights defined by \eqref{eq:fromc}. The Jacobian $D\mathcal Q(\cav)$ is similar to $W^{\mathsf T}D_{\uu'^2}$.
\item 
If $W$ is irreducible and $\uu\in\Omega_+$, the cavity solution selected in \cref{thm:selected} gives $\uu'=T_G(\uu)$ satisfying
\begin{equation}\label{eq:image}
 \rho\bigl(B(T_G(\uu)^2)\bigr)<1.
\end{equation}
\end{enumerate}
\end{proposition}
\begin{proof}
%For the general statement, 
Fix an admissible cavity solution $\cav$ and set $\beta_{v\to w}=j_{vw}/c_{v\to w}$. Equation \eqref{eq:csidentity} gives
\[
 \beta_{v\to w}=u'_{vw}\sqrt{\frac{s_w}{s_v}}.
\]
Differentiating the cavity update gives $D\mathcal Q(\cav)=W^{\mathsf T}D_{\beta^2}$. Let $L$ be diagonal on arcs with $L_{e e}=s_{o(e)}>0$, so $L$ is invertible. For an allowed transition from $f=(x\to v)$ into $e=(v\to w)$,
\[
 (L^{-1}D\mathcal Q(\cav)L)_{e f}
 =\frac1{s_v}\beta_f^2s_x=u_f'^2.
\]
Thus $L^{-1}D\mathcal Q(\cav)L=W^{\mathsf T}D_{\uu'^2}$, proving the similarity. 

Now take the solution selected in \cref{thm:selected}, for which $\uu'=T_G(\uu)$. The affine change $\zz=D_{\dd}\cav-\one$ makes the Jacobians $D\mathcal Q(\cav)$ and $D_zF_{\uu}(z_*(\uu))$ similar. By \cref{lem:iteration}, the latter has spectral radius less than one. Also, transposition and diagonal similarity give $\rho(W^{\mathsf T}D_{\uu'^2})=\rho(B(\uu'^2))$. Combining these facts with the general similarity statement yields
\[
 \rho\bigl(B(T_G(\uu)^2)\bigr)
 =\rho\bigl(D\mathcal Q(\cav)\bigr)
 =\rho\bigl(D_zF_{\uu}(z_*(\uu))\bigr)<1,
\]
which proves \eqref{eq:image}.
\end{proof}

Denoting $\cav_*$ the cavity solution selected in \cref{thm:selected},
\cref{prop:image} gives
\[
 \rho\bigl(D\mathcal Q(\cav_*)\bigr)
 =\rho\bigl(B(\uu'^2)\bigr)<1.
\]
Consequently, with the input weights $\uu$ held fixed, $\cav_*$ is a locally asymptotically stable fixed point of the iteration $\cav^{(k+1)}=\mathcal Q(\cav^{(k)})$: every initial vector sufficiently close to $\cav_*$ produces iterates converging to $\cav_*$. This stability assertion concerns the scalar cavity-coefficient update, not convergence of all messages in a Gaussian BP algorithm.

Note also that 
the condition $\rho(B(\uu'^2))<1$ involves the squared transformed weights and is distinct from the unsquared condition $\rho(B(\uu))<1$ on the original weights. The relation to Gaussian integrals is established in \cref{sec:convergence}: \cref{thm:positivity} proves that the latter condition is equivalent to positive definiteness of $H_G(\uu)$, and hence to convergence of the corresponding finite-graph Gaussian integral. Thus local stability of the cavity fixed point does not by itself guarantee convergence of that integral.

\subsection{Comparison with ratio systems}\label{sec:ratios}
To compare with \citet[Definition~1.3]{AngelFriedmanHoory2015}, put
\[
 r_{v\to w}=\frac{j_{vw}}{c_{w\to v}}.
\]
For the local operator with diagonal $h_v$ and off-diagonal $-j_{vw}$, the ratio equation at $v\to w$ reads
\[
 0=h_w-\frac{j_{vw}}{r_{v\to w}}
       -\sum_{x\sim w,\ x\ne v}j_{wx}r_{w\to x}.
\]
Substitution gives the cavity equation for $c_{w\to v}$. Each admissible solution produces a ratio system with finite, strictly positive entries $r_{v\to w}=u'_{vw}\sqrt{s_v/s_w}$. For these systems the squared-ratio matrix has entries $W_{ef}r_f^2$. For the diagonal matrix $D$ with $D_{ee}=s_{t(e)}$, it equals $D B(\uu'^2) D^{-1}$, since $t(e)=o(f)$ whenever $W_{ef}=1$. Hence the decay rate in \citet[Definition~1.3]{AngelFriedmanHoory2015} is $\alpha(r)=\rho(B(\uu'^2))$ for the ratio systems produced here.

This correspondence makes the division of labour between \cite{AngelFriedmanHoory2015} and the present results explicit. By \citet[Theorem~1.4]{AngelFriedmanHoory2015}, a pulled-back symmetric local operator on the universal covering tree has a bounded inverse exactly when it admits a ratio system with decay rate below one. Hence \cref{thm:selected,prop:image} have a consequence on the tree: for $\uu\in\Omega_+$, the selected solution furnishes a finite positive ratio system with $\alpha(r)=\rho(B(T_G(\uu)^2))<1$ for the lift of $H_G(\uu)$, so that lift is invertible on $\ell^2$ of the universal cover. In the other direction, the universal-cover theory produces, from invertibility of the lift, \emph{some} ratio system with decay rate below one; it does not produce the finite-graph statements proved here, namely positivity of the cavity coefficients and of the vertex coefficients $s_v$, the classification and uniqueness within the order region $\cav\ge\cav^{(0)}$, the monotone bracketing iteration \eqref{eq:brackets}, real-analytic dependence on $\uu$ throughout $\Omega_+$, and the congruence, prefactor, and zeta identities of \cref{thm:selected}, in which the diagonal scales are essential. Conversely, we use no classification of complex or infinite ratio systems in the proofs above.

\section{Convergence and rigidity in the positive-definite region}\label{sec:convergence}
The spectral threshold in \cref{thm:selected} involves squared weights. It must not be confused with the threshold for the ordinary weighted edge adjacency matrix $B(\uu)$. We first identify the latter with positive definiteness of the finite vertex matrix.

\subsection{Positive definiteness and the unsquared spectral radius}
The real symmetric matrix $H_G(\uu)$ has nonpositive off-diagonal entries, so it is a $Z$-matrix. In this class, positive definiteness is equivalent to being a nonsingular $M$-matrix. Such a matrix is also called a Stieltjes matrix, and its inverse is entrywise nonnegative \cite[Chapter~6]{BermanPlemmons1994}.

\begin{theorem}\label{thm:positivity}
For a finite connected simple graph and symmetric weights $0<u_e<1$,
\begin{equation}\label{eq:positivitycriterion}
 H_G(\uu)\succ0 \quad\Longleftrightarrow\quad \rho(B(\uu))<1.
\end{equation}
Irreducibility of $W$ is not required.
\end{theorem}
\begin{proof}
Suppose first that $\rho(B(\uu))<1$. For $0\le t\le1$, $B(t\uu)=tB(\uu)$, so $I-B(t\uu)$ is invertible and hence $P_G(t\uu)=\det(I-B(t\uu))\ne0$.
By \cref{prop:bass}, $H_G(t\uu)$ is nonsingular throughout this interval
since each factor $1-t^2u_e^2$ is positive.
Since $H_G(0)=I$, continuity of the eigenvalues of this real symmetric matrix implies that all its eigenvalues remain positive.

Conversely, suppose $H:=H_G(\uu)\succ0$. Its off-diagonal entries are nonpositive. 
Choose $\tau$ at least as large as the largest eigenvalue and the largest diagonal entry of $H$. Then $K=\tau I-H$ is entrywise nonnegative and $\rho(K)<\tau$. The Neumann series
\[
 H^{-1}=\tau^{-1}\sum_{r\ge0}(K/\tau)^r
\]
shows that $x=H^{-1}\one>0$. For every arc, define
\begin{equation}\label{eq:mpositivity}
 m_{v\to w}=\frac{x_v-u_{vw}x_w}{1-u_{vw}^2}.
\end{equation}
The equation $Hx=\one$, after collecting the contribution of the edge $vw$, gives
\begin{equation}\label{eq:mrecursion}
 m_{v\to w}=1+\sum_{a\sim v,\ a\ne w}u_{av}m_{a\to v}.
\end{equation}
We claim that $m>0$. If $m_{v\to w}\le0$ for some arc, then \eqref{eq:mrecursion} gives
\[
 \sum_{a\sim v,\ a\ne w}u_{av}m_{a\to v}
 =m_{v\to w}-1<0.
\]
Since every weight $u_{av}$ is positive, there is a neighbour $a\ne w$ of $v$ with $m_{a\to v}<0$. Applying the same argument to this component gives a neighbour $b\ne v$ of $a$ with $m_{b\to a}<0$. Continuing this procedure must eventually repeat an arc, since the arc set is finite. The repeated part is a directed cycle all of whose associated components are negative. For a negative component $m_{a\to v}$, \eqref{eq:mpositivity} implies $x_a<u_{av}x_v$. Multiplying these strict inequalities around the cycle gives
\[
 x_{v_0}<\left(\prod_{e\text{ on the cycle}}u_e\right)x_{v_0}<x_{v_0},
\]
a contradiction. Thus $m>0$.

Equation \eqref{eq:mrecursion} now reads
\[
 W^{\mathsf T}D_{\uu}m=m-\one<m.
\]
Conjugation by $\diag(m)$ bounds the spectral radius by the maximum of these component ratios, which is strictly less than one. Finally $\rho(W^{\mathsf T}D_{\uu})=\rho(WD_{\uu})$, by transposition and diagonal similarity. This proves the reverse implication.
\end{proof}
The forward implication is a specialization of \citet[Theorem~10]{WatanabeFukumizu2011}. Their spectral hypothesis is that the directed-edge matrix has no eigenvalue in $\RR_{\ge1}$, a condition that is in general weaker than requiring its spectral radius to be less than one. In the present positive-weight setting, however, $B(\uu)$ is nonnegative, so the Perron--Frobenius theorem ensures that $\rho(B(\uu))$ itself is a nonnegative real eigenvalue, without an irreducibility assumption. Consequently,
the condition $\Spec(B(\uu))\cap\RR_{\ge1}=\varnothing$
is equivalent to $\rho(B(\uu))<1$.
Thus the two spectral conditions are equivalent for the weights considered here. The homotopy argument above follows their proof. The additional implication established in \cref{thm:positivity} is the reverse direction, which uses the positivity, symmetry under reversal, and strict upper bound $u_e<1$ of the present weights; it is not asserted for general complex or signed weights.

\subsection{Domains of convergence for the Euler product and Gaussian integral representations}
For complex weights, a sufficient condition for absolute convergence of the logarithmic Euler expansion is
\begin{equation}\label{eq:eulerconvergence}
 \rho\bigl(B(|\uu|)\bigr)<1.
\end{equation}
Indeed, the absolute values of cycle contributions are bounded by those for the nonnegative matrix $B(|\uu|)$, and $\sum_{k\ge1}\Tr(B(|\uu|)^k)/k$ converges under \eqref{eq:eulerconvergence}. Two convenient sufficient conditions are
\begin{equation}\label{eq:rowcondition}
 \max_{v\to w}\sum_{x\sim w,\ x\ne v}|u_{wx}|<1,
 \qquad\text{or, more restrictively,}\qquad
 (d_{\max}-1)\max_e|u_e|<1,
\end{equation}
where $d_{\max}:=\max_{v\in V}d_v$ is the maximum vertex degree of $G$.
They follow from row-sum bounds on $B(|\uu|)$.

For real positive weights, the Gaussian integral
\begin{equation}\label{eq:GaussianIntegral}
 \mathcal Z_G(\uu)=
 \frac{\prod_{e\in E}(1-u_e^2)^{-1}}{\pi^n}
 \int_{\CC^n}e^{-X^\dagger H_G(\uu)X}\prod_{v\in V}d^2X_v
\end{equation}
is finite exactly when $H_G(\uu)\succ0$. Diagonalizing the real symmetric matrix, and integrating the real and imaginary parts independently, then gives
\begin{equation}\label{eq:GaussianEqualsZeta}
 \mathcal Z_G(\uu)=\frac{\prod_e(1-u_e^2)^{-1}}{\det H_G(\uu)}
 =\zeta_G(\uu).
\end{equation}
If an eigenvalue is zero or negative, the integral diverges. In this case, the integral representation in \eqref{eq:GaussianEqualsZeta} does not apply, and we do not assign a finite value to the divergent integral by regularization.

An alternative local sufficient condition for positivity is
\begin{equation}\label{eq:diagdom}
 \sum_{w\sim v}\frac{u_{vw}}{1+u_{vw}}<1
 \quad\text{for every }v.
\end{equation}
In fact,
\[
 H_{vv}-\sum_{w\ne v}|H_{vw}|
 =1-\sum_{w\sim v}\frac{u_{vw}}{1+u_{vw}}.
\]
Thus \eqref{eq:diagdom} makes the real symmetric matrix strictly diagonally dominant with positive diagonal. The criteria \eqref{eq:rowcondition} and \eqref{eq:diagdom} need not be sharp; \cref{thm:positivity} supplies the exact spectral condition in this positive setting.

\subsection{Rigidity of positive diagonal congruence}
\begin{theorem}[Rigidity]\label{thm:rigidity}
Let $\uu,\uu'\in(0,1)^E$ and $S=\diag(a_v)>0$ satisfy
\[
 S H_G(\uu')S=H_G(\uu).
\]
If $H_G(\uu)\succ0$, then $S=I$ and $\uu'=\uu$.
\end{theorem}
\begin{proof}
By \cref{thm:correspondence}, the reference solution
\[
 c^{(0)}_{v\to w}=\frac1{1-u_{vw}^2}
\]
and the solution
\[
 c^{(1)}_{v\to w}=\frac{a_v^2}{1-u_{vw}'^{\,2}}
\]
satisfy the same cavity equations. Subtract those equations and set $d=c^{(1)}-c^{(0)}$. Then
\begin{equation}\label{eq:cavitydifference}
 d_{v\to w}=\sum_{x\sim v,\ x\ne w}
 \frac{j_{xv}^2}{c^{(0)}_{x\to v}c^{(1)}_{x\to v}}d_{x\to v}.
\end{equation}
The off-diagonal entries of the congruence imply
\[
 a_xa_v\frac{u'_{xv}}{1-u_{xv}'^{\,2}}
 =\frac{u_{xv}}{1-u_{xv}^2}.
\]
Substitution into \eqref{eq:cavitydifference} gives
\[
 \frac{j_{xv}^2}{c^{(0)}_{x\to v}c^{(1)}_{x\to v}}
 =\frac{a_v}{a_x}u_{xv}u'_{xv}.
\]
Set $A:=W^{\mathsf T}D_{\uu\uu'}$ and $\gamma:=\max_e u'_e<1$. Since $u_eu'_e\le\gamma u_e$ for every edge and $W$ is nonnegative, we have the entrywise comparison
\[
 0\le A\le\gamma W^{\mathsf T}D_{\uu}.
\]
By the assumption $H_G(\uu)\succ0$ and \cref{thm:positivity}, $\rho(B(\uu))<1$. Monotonicity of the spectral radius for nonnegative matrices therefore gives
\[
 \rho(A)\le\gamma\rho(W^{\mathsf T}D_{\uu})
 =\gamma\rho(B(\uu))<1.
\]
%Thus $A$ has no eigenvalue equal to one: a nonzero vector satisfying $f=Af$ would be an eigenvector with eigenvalue one and would contradict $\rho(A)<1$.
Now define $f_{v\to w}:=d_{v\to w}/a_v$. The difference equation above then becomes
\begin{equation}\label{eq:differencespectral}
 f=W^{\mathsf T}D_{\uu\uu'}f.
\end{equation}
Hence $f=0$. Since every $a_v>0$, this implies $d=0$ and therefore $c^{(1)}=c^{(0)}$. Their vertex coefficients agree: $a_v^2=s_v(c^{(1)})=s_v(c^{(0)})=1$. The edge congruence and the strict monotonicity of $u/(1-u^2)$ on $(0,1)$ now imply $u'_e=u_e$.
\end{proof}

The assumption of positive edge weights in \cref{thm:rigidity} is essential: allowing negative edge weights permits nontrivial positive diagonal congruences even when both weighted vertex matrices are positive definite. Thus positive definiteness alone does not force the weight transformation to be the identity. See \cref{sec:signedexample} for an explicit example and its implications for convergence.

\begin{corollary}\label{cor:noconvergence}
For $\uu\in\Omega_+$, neither $H_G(\uu)$ nor $H_G(T_G(\uu))$ is positive definite. Consequently the nontrivial selected identity is not an equality of two convergent positive Gaussian integrals of the form \eqref{eq:GaussianIntegral}.
\end{corollary}
\begin{proof}
Let $\uu':=T_G(\uu)$. By \cref{thm:selected}, $0<u'_e<u_e<1$ for every edge, so this weight transformation is nontrivial. If $H_G(\uu)$ were positive definite, \cref{thm:rigidity} would force $\uu'=\uu$ and $S_G(\uu)=I$, which gives a contradiction. Thus the original matrix $H_G(\uu)$ is not positive definite.

For the transformed matrix, we use the congruence from \eqref{eq:selectedcongruence},
\[
 S_G(\uu)H_G(\uu')S_G(\uu)=H_G(\uu).
\]
The matrix $S_G(\uu)$ is invertible, so the two matrices have the same inertia, namely the same numbers of positive, negative, and zero eigenvalues. Since $H_G(\uu)$ is not positive definite, $H_G(\uu')$ cannot be positive definite either. 
%By the convergence criterion for \eqref{eq:GaussianIntegral}, both corresponding Gaussian integrals therefore diverge.
\end{proof}
Note that the failure of positive definiteness for the original matrix also follows directly from \cref{thm:positivity}. 
In fact, 
since $0<u_e^2<u_e$, monotonicity of the spectral radius gives
\[
 \rho(B(\uu))\ge\rho(B(\uu^2))>1.
\]
Therefore, $H_G(\uu)$ cannot be positive definite from \cref{thm:positivity}.
For the transformed weights $\uu'=T_G(\uu)$, however, the condition $\rho(B(\uu'^2))<1$ from \cref{prop:image} alone does not determine whether $\rho(B(\uu'))<1$. The nonpositive-definiteness of $H_G(\uu')$ follows from the congruence in the proof, rather than from the squared-weight inequality alone.

This corollary does not invalidate the determinant identity. It identifies the correct interpretation of that identity. Positive inertia, negative inertia, and nullity are preserved by every congruence, and the positive prefactor in \eqref{eq:prefactorintro} preserves the zero set of $P_G$ along the specified transformation.

\begin{figure}[t]
\centering
\begin{tikzpicture}[x=11.5cm,y=1cm]
  \fill[black!5]  (0,0) rectangle (0.5,0.5);
  \fill[black!14] (0.5,0) rectangle (0.70711,0.5);
  \fill[black!26] (0.70711,0) rectangle (1,0.5);
  \draw (0,0) rectangle (1,0.5);
  \draw (0.5,0) -- (0.5,0.5);
  \draw (0.70711,0) -- (0.70711,0.5);
  \node[below] at (0,-0.05) {$0$};
  \node[below] at (0.5,-0.05) {$1/t$};
  \node[below] at (0.70711,-0.05) {$1/\sqrt{t}$};
  \node[below] at (1,-0.05) {$1$};
  \node[right] at (1.02,0.25) {$u$};
  \node[align=center,font=\small] at (0.25,1.45)
    {$H_G(u\one)\succ0$\\ $\rho(B(u\one))<1$\\ rigidity (\cref{thm:rigidity})};
  \node[align=center,font=\small] at (0.6036,1.45)
    {$\rho(B(u\one))\ge1$\\ $\rho(B(u^2\one))\le1$\\ contains $u'$};
  \node[align=center,font=\small] at (0.8536,1.45)
    {$\Omega_+$: $\rho(B(u^2\one))>1$\\ domain of $T_G$};
  \draw[dotted] (0.9,0) -- (0.9,-0.75);
  \draw[dotted] (0.5556,0) -- (0.5556,-0.75);
  \draw[->,thick] (0.9,-0.75) .. controls (0.82,-1.35) and (0.64,-1.35) .. (0.5556,-0.75);
  \node[below,font=\small] at (0.728,-1.2) {$T_G:\ u\mapsto u'=1/(tu)$};
\end{tikzpicture}
\caption{Spectral regions for uniform weights $u_e=u$ on a $(t+1)$-regular graph with $t>1$, drawn for $t=2$. Here $\rho(B(u\one))=tu$ and $\rho(B(u^2\one))=tu^2$. For $1/\sqrt t<u<1$ the selected transformation is $u'=1/(tu)$, so $1/t<u'<1/\sqrt t$. The reciprocal formula extends to the fixed point $1/\sqrt t$ on the boundary of the domain (\cref{rem:boundary}). For independent edge weights the conclusion is the inclusion in \cref{rem:band}.}
\label{fig:regions}
\end{figure}
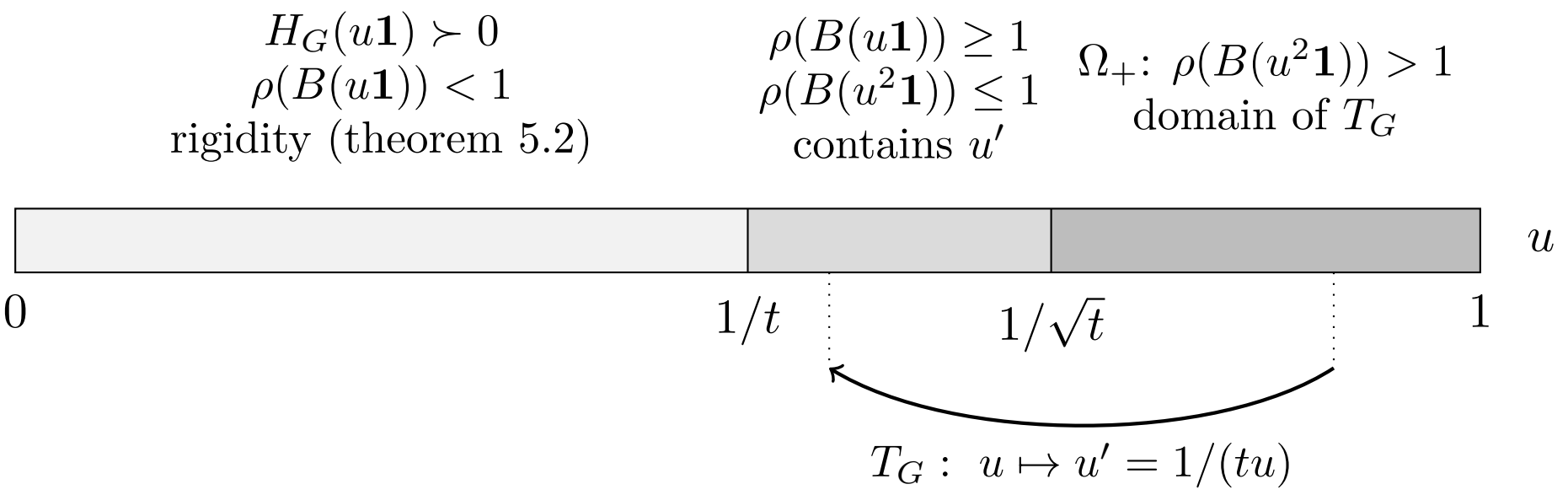

\begin{remark}[The intermediate region]\label{rem:band}
\Cref{thm:positivity,prop:image,cor:noconvergence} give the inclusion
\[
 T_G(\Omega_+)\ \subseteq\
 \bigl\{\vv\in(0,1)^E:\ \rho(B(\vv^2))<1\le\rho(B(\vv))\bigr\}.
\]
Thus every output has subcritical squared weights and a vertex matrix that fails to be positive definite. \Cref{fig:regions} illustrates the inclusion for uniform weights on a regular graph.
\end{remark}

\section{Uniform weights and vertex-factorized subfamilies}\label{sec:examples}
\subsection{A complete classification for uniform input and output}
The classical regular and biregular reciprocal relations are recovered by restricting both input and output to uniform weights. This also shows exactly what is lost when one insists on only one parameter.

\begin{proposition}\label{prop:uniform}
Suppose $u,u'\in(0,1)$, $u'\ne u$, and $u_e=u$, $u'_e=u'$ for every edge of a connected graph; uniform weight vectors are abbreviated by their common value, so that for instance $\zeta_G(u)$ stands for $\zeta_G(u\one)$. A positive diagonal congruence exists if and only if
\begin{equation}\label{eq:degreeproduct}
 (d_v-1)(d_w-1)=Q:=\frac1{u^2u'^2}
 \quad\text{for every }vw\in E.
\end{equation}
It follows that the graph is either regular or bipartite with constant degree on each bipartition class. Conversely, these graphs give the required congruence whenever $u'=1/(\sqrt Q\,u)$ and $u,u'\in(0,1)$. The diagonal scales are uniquely determined by
\begin{equation}\label{eq:uniformscale}
 a_v^2=\frac{1-u'^2}{1-u^2}
 \frac{1+(d_v-1)u^2}{1+(d_v-1)u'^2}.
\end{equation}
\end{proposition}
\begin{proof}
Matching an edge entry and a diagonal entry gives, respectively,
\[
 a_va_w=\frac u{u'}\frac{1-u'^2}{1-u^2},\qquad
 a_v^2=\frac{1-u'^2}{1-u^2}
 \frac{1+(d_v-1)u^2}{1+(d_v-1)u'^2}.
\]
Multiply the diagonal equations at the ends of an edge and compare with the square of the edge equation. If $t_v=d_v-1$, the result is
\[
 u'^2(1+t_vu^2)(1+t_wu^2)
 -u^2(1+t_vu'^2)(1+t_wu'^2)
 =(u'^2-u^2)(1-t_vt_wu^2u'^2)=0.
\]
Since $u'\ne u$, this is \eqref{eq:degreeproduct}. Conversely, \eqref{eq:degreeproduct} and the positive square roots in \eqref{eq:uniformscale} satisfy both edge and diagonal equations.

All $t_v$ are positive. Along a path, $t_v$ alternates between a value $r$ and $Q/r$. An odd simple cycle forces $r^2=Q$ and hence regularity by connectedness. If there is no odd simple cycle, the graph is bipartite and $t_v$ is constant on each class. This proves the classification.
\end{proof}

This proposition classifies positive diagonal congruences with uniform input and output, not all possible functional equations of a graph zeta function. The regular reciprocal relation is classical \cite{Terras2011}. The biregular relation belongs to Hashimoto's theory \cite{Hashimoto1989} and is applied to hypergraph zeta functions in \cite{Storm2006}.

For a $(t+1)$-regular graph,
\begin{equation}\label{eq:regularrelation}
 u'=\frac1{tu},\qquad
 a^2=\frac{t^2u^2-1}{t(1-u^2)},\qquad
 \zeta_G(u')=(tu^2)^m
 \left(\frac{t(1-u^2)}{t^2u^2-1}\right)^{m-n}\zeta_G(u).
\end{equation}
The real positive scales require $u>1/t$. The nonnegative branch of \cref{thm:selected} is the part with $u>1/\sqrt t$, for which $u'<u$. The reciprocal direction is represented by a different cavity order region. The equal-weight reciprocal map, continued algebraically, is an involution; the selected many-weight map is not asserted to be one.

If the two bipartite classes have degrees $r+1$ and $s+1$, then $Q=rs$ and
\begin{align}
 u'&=\frac1{\sqrt{rs}\,u},\label{eq:biregularp}\\
 a_A^2&=\frac{rsu^2-1}{r(1-u^2)}\frac{1+ru^2}{1+su^2},&
 a_B^2&=\frac{rsu^2-1}{s(1-u^2)}\frac{1+su^2}{1+ru^2}.
 \label{eq:biregularscales}
\end{align}
The prefactor follows from \eqref{eq:prefactorintro} without any extra choice of normalization.

\begin{example}[$K_{2,3}$]\label[example]{ex:k23}
The classes have degrees $3$ and $2$. For $u^2=3/4$ and $u'^2=2/3$,
\[
 a_A^2=\frac{10}{7},\qquad a_B^2=\frac75,\qquad
 C_G=\frac{5103}{5120}.
\]
An independent determinant calculation gives
\[
 P_{K_{2,3}}(u)=(1-u^4)^2(1-4u^4),
\]
which verifies the relation. The exact calculation is described in \cref{sec:exampledata}. Here the two scales differ, although every edge has the same input and output weight.
\end{example}

\subsubsection{The diamond graph: uniform input and nonuniform output}\label{sec:diamondexample}
\begin{figure}[t]
\centering
\begin{tikzpicture}[scale=1.05,vtx/.style={circle,fill,inner sep=1.5pt},lab/.style={font=\scriptsize}]
\begin{scope}
 \node[vtx,label={[font=\small]left:$1$}]  (a1) at (0,0) {};
 \node[vtx,label={[font=\small]above:$2$}] (b1) at (1.4,1.1) {};
 \node[vtx,label={[font=\small]right:$3$}] (c1) at (2.8,0) {};
 \node[vtx,label={[font=\small]below:$4$}] (d1) at (1.4,-1.1) {};
 \draw (a1)--(b1) node[midway,above left,lab]{$0.85$};
 \draw (b1)--(c1) node[midway,above right,lab]{$0.85$};
 \draw (c1)--(d1) node[midway,below right,lab]{$0.85$};
 \draw (d1)--(a1) node[midway,below left,lab]{$0.85$};
 \draw (a1)--(c1) node[midway,above,lab]{$0.85$};
\end{scope}
\draw[->,very thick] (3.5,0) -- (4.7,0) node[midway,above,font=\small]{$T_G$};
\begin{scope}[shift={(5.4,0)}]
 \node[vtx,label={[font=\small]left:$1$}]  (a2) at (0,0) {};
 \node[vtx,label={[font=\small]above:$2$}] (b2) at (1.4,1.1) {};
 \node[vtx,label={[font=\small]right:$3$}] (c2) at (2.8,0) {};
 \node[vtx,label={[font=\small]below:$4$}] (d2) at (1.4,-1.1) {};
 \draw (a2)--(b2) node[midway,above left,lab]{$0.7741$};
 \draw (b2)--(c2) node[midway,above right,lab]{$0.7741$};
 \draw (c2)--(d2) node[midway,below right,lab]{$0.7741$};
 \draw (d2)--(a2) node[midway,below left,lab]{$0.7741$};
 \draw (a2)--(c2) node[midway,above,lab]{$0.7704$};
\end{scope}
\end{tikzpicture}
\caption{The worked functional equation of \cref{ex:diamond} on the diamond graph. Left: uniform input weight $0.85$ on every edge. Right: the output weights $T_G(\uu)$, rounded to four digits, equal on the four outer edges and different on the diagonal. The accompanying vertex scales are $a_1=a_3=1.271404$ and $a_2=a_4=1.247134$, and $\zeta_G(\uu')=C_G\,\zeta_G(\uu)$ with $C_G=0.992174223093$.}
\label{fig:diamond}
\end{figure}
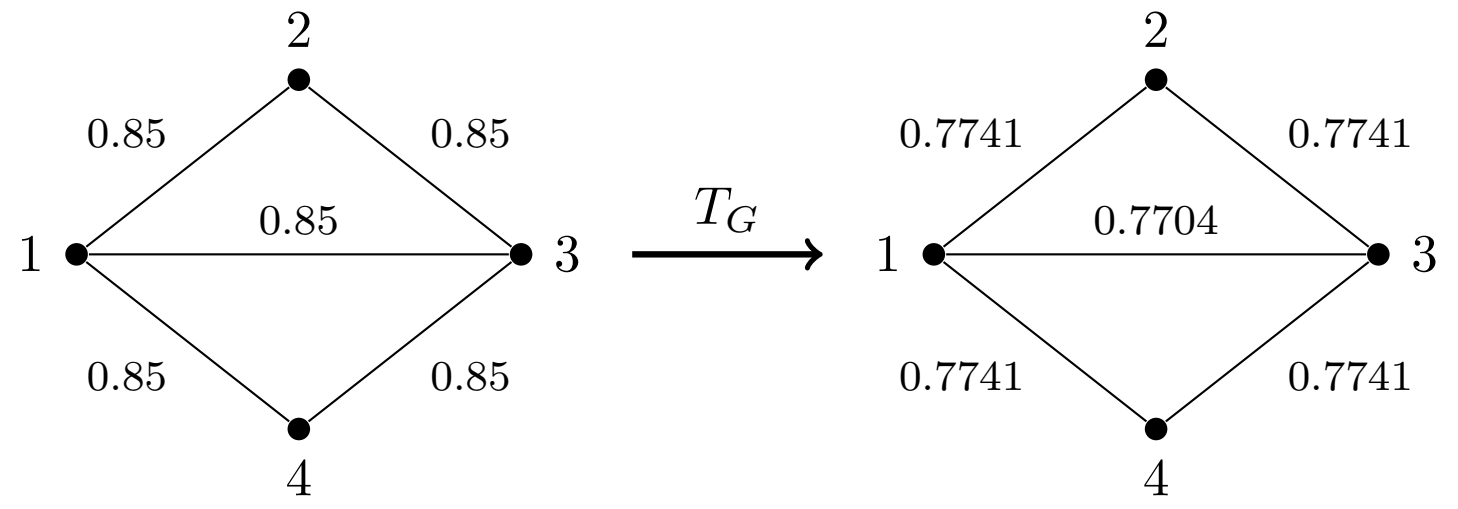

\begin{example}\label[example]{ex:diamond}
Let $G$ have edges $12,23,34,41,13$, so that vertices $1,3$ have degree $3$ and vertices $2,4$ degree $2$ (\cref{fig:diamond}). The product $(d_v-1)(d_w-1)$ is $2$ on an outer edge and $4$ on the diagonal $13$. Thus \cref{prop:uniform} rules out a nontrivial positive diagonal congruence with both input and output uniform on this graph.

The weight transformation of \cref{thm:main} gives a nonuniform output in this example. Here $\rho(W)=1.521380\ldots$, so the uniform input $u_e=0.85$ lies in $\Omega_+$, with $\rho(B(\uu^2))=0.85^2\rho(W)=1.099197\ldots$. The monotone iteration in \cref{sec:selectedalgorithm}, with the input data in \cref{sec:exampledata}, solves \eqref{eq:mainsystem} and gives the nonuniform output weights and vertex scales below. The symmetries $1\leftrightarrow3$ and $2\leftrightarrow4$ also allow a reduction to three unknowns, although the algorithm uses all directed edges.
\[
 u'_{12}=u'_{23}=u'_{34}=u'_{41}=0.774125257693,\qquad
 u'_{13}=0.770362021074,
\]
\[
 s_1=s_3=1.616467117614,\qquad s_2=s_4=1.555343880069,\qquad
 C_G=0.992174223093,
\]
every displayed digit being rounded. The resulting instance of the functional equation \eqref{eq:mainintro} reads, numerically,
\[
 \underbrace{\zeta_G(\uu')}_{-4.417257427525}
 \;=\;\underbrace{C_G}_{0.992174223093}\;\times\;\underbrace{\zeta_G(\uu)}_{-4.452098557605},
\]
equivalently $P_G(\uu)=-0.224613176699$ and $P_G(\uu')=-0.226384813746$ with $P_G(\uu')\,C_G=P_G(\uu)$ to machine precision. Both zeta values are negative: the weights lie outside the positive-definite region, consistently with \cref{cor:noconvergence}, and the identity is one of rational functions. The output satisfies $\rho(B(\uu'^2))=0.909756\ldots<1$, as \cref{prop:image} requires. Thus failure of the one-parameter ansatz does not prevent the many-weight construction. The independent residual checks are specified in \cref{sec:checkalgorithm}.
\end{example}

\subsection{Weights initially assigned to vertices}\label{sec:vertexinputs}
A natural vertex-weighted specialization is
\begin{equation}\label{eq:vertexweights}
 u_{vw}=\sqrt{q_vq_w},\qquad q_v>0.
\end{equation}
For a cycle, $\prod_eu_e$ equals the product of $q_v$ over the visited vertices, counted with multiplicity. Vertex-weighted Ihara and Bartholdi variants have been studied in other conventions by \cite{KonnoEtAl2019,Zhu2021}; \eqref{eq:vertexweights} specifies the convention used here. The letter $q$ is reserved throughout for such vertex factors; edge weights are always denoted $u$.

The many-weight transformation need not be restricted in advance to this family. There is a precise test for whether a given output belongs to it. Let $L$ be the unsigned edge--vertex incidence matrix, with $L_{e v}=1$ if $v$ is an endpoint of $e$ and $0$ otherwise. Writing $b_v=\log q'_v$, factorization is equivalent to
\begin{equation}\label{eq:vertexclosure}
 Lb=2\log\uu'.
\end{equation}
Thus positive vertex factors exist precisely when $2\log\uu'\in\im L$. To require $0<q'_v<1$ as well, one must in addition find a solution with every $b_v<0$. The inequality $u'_e<1$ alone does not imply this latter requirement.

For example, on $K_4$ equation \eqref{eq:vertexclosure} requires
\begin{equation}\label{eq:K4closure}
 u'_{12}u'_{34}=u'_{13}u'_{24}=u'_{14}u'_{23}.
\end{equation}
These conditions are also sufficient for positive vertex factors on $K_4$, as follows by solving for $b_1,b_2,b_3$ from a triangle and using \eqref{eq:K4closure} for the remaining edges. The example in \cref{sec:inhomogeneous} illustrates that the selected cavity output need not satisfy these relations. The transformation theorem itself always concerns independent edge weights.

\subsubsection{\texorpdfstring{$K_4$}{K4}: loss of vertex factorization}\label{sec:inhomogeneous}
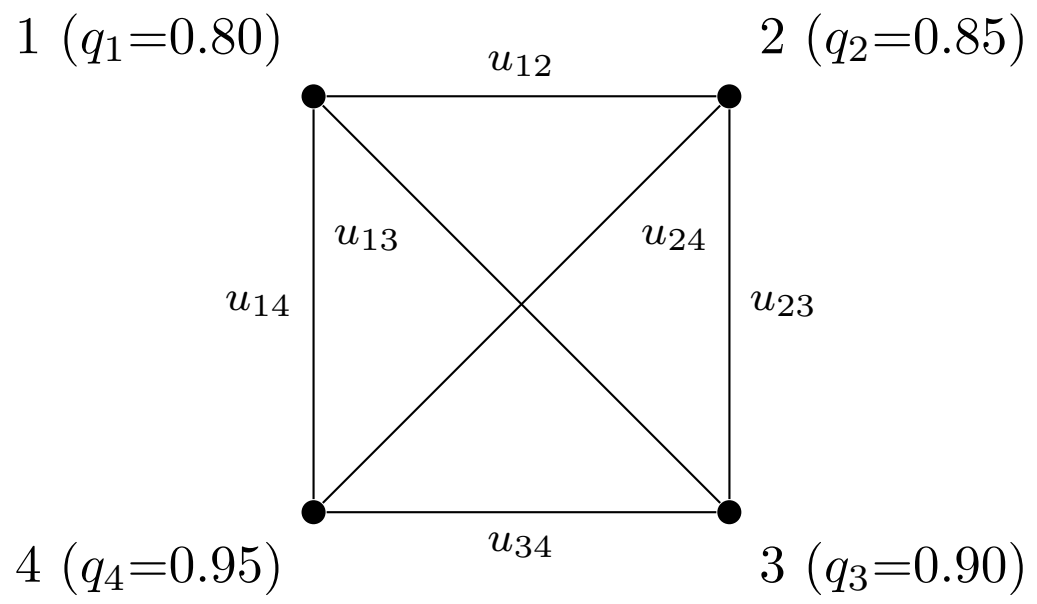
\begin{figure}[t]
\centering
\begin{tikzpicture}[scale=1.25,vtx/.style={circle,fill,inner sep=1.6pt},lab/.style={font=\scriptsize}]
 \node[vtx,label={[font=\small]above left:{$1\ (q_1{=}0.80)$}}]  (v1) at (0,2.2) {};
 \node[vtx,label={[font=\small]above right:{$2\ (q_2{=}0.85)$}}] (v2) at (2.2,2.2) {};
 \node[vtx,label={[font=\small]below right:{$3\ (q_3{=}0.90)$}}] (v3) at (2.2,0) {};
 \node[vtx,label={[font=\small]below left:{$4\ (q_4{=}0.95)$}}]  (v4) at (0,0) {};
 \draw (v1)--(v2) node[midway,above,lab]{$u_{12}$};
 \draw (v2)--(v3) node[midway,right,lab]{$u_{23}$};
 \draw (v3)--(v4) node[midway,below,lab]{$u_{34}$};
 \draw (v4)--(v1) node[midway,left,lab]{$u_{14}$};
 \draw (v1)--(v3) node[pos=0.25,below left,lab]{$u_{13}$};
 \draw (v2)--(v4) node[pos=0.25,below right,lab]{$u_{24}$};
\end{tikzpicture}
\caption{The complete graph $K_4$ used in \cref{sec:inhomogeneous}, with the vertex factors $q_i$ that generate the input edge weights $u_{ij}=\sqrt{q_iq_j}$ of \eqref{eq:vertexweights}. The numerical values of the $u_{ij}$ are listed in the text; the transformed weights $u'_{ij}=T_G(\uu)_{ij}$ appear in \cref{tab:inhomogeneous}.}
\label{fig:k4example}
\end{figure}

This example shows that the selected transformation need not preserve vertex factorization. The graph is the complete graph $K_4$ on the vertex set $\{1,2,3,4\}$, with all six edges present; see \cref{fig:k4example}. The input weights are generated from vertex factors as in \eqref{eq:vertexweights}: set
\[
 (q_1,q_2,q_3,q_4)=(0.80,0.85,0.90,0.95),\qquad
 u_{ij}=\sqrt{q_iq_j},
\]
so that, rounded to six digits,
\begin{align*}
 u_{12}&=0.824621, & u_{13}&=0.848528, & u_{14}&=0.871780,\\
 u_{23}&=0.874643, & u_{24}&=0.898610, & u_{34}&=0.924662.
\end{align*}
We used $80$-digit arithmetic for the fixed point and double precision for the spectral radii, as specified in \cref{sec:computations,sec:selectedalgorithm,sec:checkalgorithm}. The rounded values below are numerical checks at those precisions. The input lies in $\Omega_+$, since $\rho(B(\uu^2))=1.527083\ldots>1$. Iterating \eqref{eq:eta} from an upper bound, as specified in \cref{sec:selectedalgorithm}, converges to the selected branch of \cref{thm:main}. The six output weights $u'_{ij}=T_G(\uu)_{ij}$ are listed in \cref{tab:inhomogeneous}; the vertex scales are
\[
 (s_1,s_2,s_3,s_4)=(4.224295,\ 4.398475,\ 4.523994,\ 4.597532),
 \qquad a_v=\sqrt{s_v},
\]
and the prefactor is $C_G=0.668069764833$. The identity verified is the functional equation \eqref{eq:mainintro} for these data: evaluating the two sides independently gives
\[
 P_G(\uu)=-0.198770933196,\qquad
 P_G(\uu')=-0.297530203669,
\]
with $\bigl|P_G(\uu')\,C_G-P_G(\uu)\bigr|<10^{-15}$; the congruence $SH_G(\uu')S=H_G(\uu)$ with $S=\diag(\sqrt{s_v})$ holds entrywise with residual below $10^{-14}$. The output satisfies $\rho(B(\uu'^2))=0.657857\ldots<1$, as \cref{prop:image} requires.

\begin{table}[htbp]
\centering
\caption{Selected output for a vertex-factorized input on $K_4$. Entries are rounded numerical values, not exact identities.}\label{tab:inhomogeneous}
\begin{tabular}{cc@{\qquad}cc}
\toprule
Edge & $u'_{ij}$ & Edge & $u'_{ij}$\\
\midrule
$12$ & $0.467287648309$ & $23$ & $0.566547203052$\\
$13$ & $0.511704825468$ & $24$ & $0.628291857183$\\
$14$ & $0.562991973839$ & $34$ & $0.704426381885$\\
\bottomrule
\end{tabular}
\end{table}
The products entering \eqref{eq:K4closure} are approximately
\[
 0.329169747398,\qquad 0.321499975122,\qquad 0.318961528119,
\]
respectively. Since they differ, the output is not vertex-factorized although the input is; the exact criterion remains \eqref{eq:vertexclosure}. The procedure in \cref{sec:checkalgorithm} checks the cavity equations, the diagonal congruence, both determinant representations, and \eqref{eq:prefactorintro}, with the two sides evaluated independently.

\section{Exponentially many positive cavity branches}\label{sec:manybranches}
The uniqueness assertion in \cref{thm:main} concerns only solutions satisfying
\[
 z_e=(1-u_e^2)c_e-1>0,
 \qquad\text{equivalently}\qquad
 c_e>c_e^{(0)}=\frac{1}{1-u_e^2}
 \quad(e\in\arcs).
\]
In this section we require only admissibility in the sense of \cref{def:admissible}: $c_e>0$ for every arc and $s_v>0$ for every vertex. These positivity conditions do not require $\cav\ge\cav^{(0)}$.

Under the hypotheses of \cref{thm:main}, \cref{thm:selected} shows that the only cavity solutions in $\cav\ge\cav^{(0)}$ are the reference solution and the unique selected solution. Every additional admissible solution must therefore have $c_e<c_e^{(0)}$, equivalently $z_e<0$, on at least one arc. The construction below produces many such additional solutions, consistently with the uniqueness in \cref{thm:main}.

We first prove the construction in \cref{sec:branchconstruction}, then illustrate its branches for $k=2$ in \cref{sec:fourbranches}.

\subsection{Construction and existence of the branches}\label{sec:branchconstruction}
\begin{theorem}\label{thm:manybranches}
For every integer $k\ge1$, there is a connected simple graph $G_k$ with positive weights less than one for which the cavity equations have at least $2^k$ distinct solutions with all $c_e>0$ and all $s_v>0$. Hence $G_k$ admits at least $2^k$ distinct pairs $(\uu',S)$ in \cref{thm:correspondence}. The count is for labelled graphs and does not identify solutions under graph automorphisms.
\end{theorem}

\begin{figure}[t]
\centering
\begin{tikzpicture}[scale=0.85,vtx/.style={circle,fill,inner sep=1.5pt}]
  \foreach \k in {0,1,2}{
    \begin{scope}[shift={(4.2*\k,0)}]
      \node[vtx] (a\k) at (0,0) {};
      \node[vtx] (b\k) at (1.1,1) {};
      \node[vtx] (c\k) at (1.1,-1) {};
      \node[vtx] (d\k) at (2.2,0) {};
      \draw (a\k)--(b\k) (a\k)--(c\k) (a\k)--(d\k) (b\k)--(c\k) (b\k)--(d\k) (c\k)--(d\k);
    \end{scope}
  }
  \draw[dashed] (d0)--(a1) node[midway,above] {$\varepsilon$};
  \draw[dashed] (d1)--(a2) node[midway,above] {$\varepsilon$};
\end{tikzpicture}
\caption{The chain $G_k$ of \cref{thm:manybranches}, drawn for $k=3$: labelled copies of $K_4$ with weight $4/5$ on every solid edge, joined by bridges of weight $\varepsilon$ (dashed). At $\varepsilon=0$ each copy independently selects one of the two constant cavity roots \eqref{eq:K4roots}, and each of the resulting $2^k$ nondegenerate solutions continues to small $\varepsilon>0$.}
\label{fig:chain}
\end{figure}
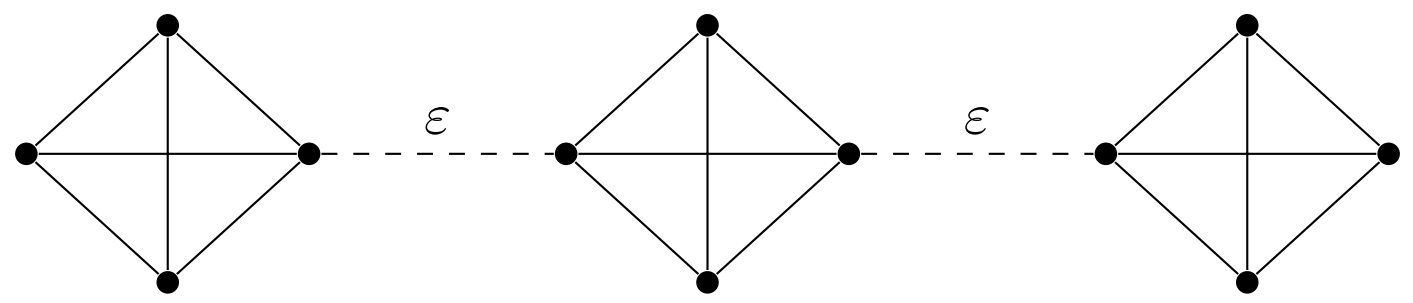

\begin{proof}
Take $k$ disjoint labelled copies of $K_4$ and connect successive copies by one bridge, obtaining a chain (\cref{fig:chain}). Every internal edge carries the weight $4/5$, and every bridge the weight $\varepsilon>0$. Let $h_v(\varepsilon)$ and $j_{vw}(\varepsilon)$ be the coefficients in \eqref{eq:hjintro} for these weights, and let $\mathcal Q_\varepsilon$ be the corresponding cavity update. The cavity equations \eqref{eq:cavityupdate}, imposed on all arcs of $G_k$, including the bridge arcs, are
\[
 F_{v\to w}(\cav,\varepsilon)
 :=c_{v\to w}-h_v(\varepsilon)
 +\sum_{x\sim v,\ x\ne w}
 \frac{j_{xv}(\varepsilon)^2}{c_{x\to v}}=0
 \qquad(v\to w\in\arcs).
\]
Thus $F(\cav,\varepsilon):=(F_{v\to w}(\cav,\varepsilon))_{v\to w\in\arcs}=\cav-\mathcal Q_\varepsilon(\cav)$ is the residual map of the cavity equations, with unknown $\cav$ and parameter $\varepsilon$. For the existence argument, extend these equations analytically to $\varepsilon=0$. This is an auxiliary boundary point; the final weights will be strictly positive.

Inside each copy of $K_4$, at $\varepsilon=0$,
\[
 h=\frac{19}{3},\qquad j=\frac{20}{9}.
\]
A constant cavity value satisfies $c=h-2j^2/c$, whose two roots are
\begin{equation}\label{eq:K4roots}
 c_0=\frac{25}{9},\qquad c_1=\frac{32}{9}.
\end{equation}
Their full vertex coefficients and output weights are
\begin{equation}\label{eq:K4rootdata}
 s_0=1,\quad s_1=\frac{13}{6},\qquad
 u'_0=\frac45,\quad u'_1=\frac58.
\end{equation}
Each copy independently chooses $c_0$ or $c_1$. An outgoing bridge cavity coefficient at $\varepsilon=0$ equals $s_0$ or $s_1$ of its source copy. Incoming bridge contributions vanish because the bridge coupling $j_{\mathrm{bridge}}$ is zero. Hence there are $2^k$ distinct positive solutions at this boundary point, all with positive $s_v$.

Before applying the implicit-function theorem,
we check that each of these $2^k$ solutions at $\varepsilon=0$ is a nondegenerate solution of the cavity equations: the Jacobian $D_{\cav}F(\cav,0)=I-D_{\cav}\mathcal Q_0(\cav)$ must be invertible at each such solution.

Fix a choice $(r_1,\ldots,r_k)\in\{0,1\}^k$, where copy $i$ selects the root $c_{r_i}$. Order the internal arc variables by copy, followed by all bridge arc variables. All Jacobians in the following calculation are evaluated at the corresponding solution with $\varepsilon=0$. Since bridge terms carry $j_{\mathrm{bridge}}^2=0$, every column corresponding to a bridge variable vanishes. Thus
\[
 D\mathcal Q=
 \begin{pmatrix}A&0\\ C&0\end{pmatrix},
 \qquad
 I-D\mathcal Q=
 \begin{pmatrix}I-A&0\\ -C&I\end{pmatrix},
\]
where $A$ and $C$ describe the dependence of internal updates on internal variables and the dependence of bridge updates on internal variables, respectively. The block $C$ need not vanish; the identity block in $I-D\mathcal Q$ acts on the bridge variables. Different copies have no direct interaction in the internal updates at $\varepsilon=0$, so $A$ is block diagonal, with one $12\times12$ block for each copy:
\[
 \begin{gathered}
 A=\diag(A_1,\ldots,A_k),\qquad
 A_i=\alpha_{r_i}W_{K_4}^{\mathsf T},\\
 \alpha_0=\frac{j^2}{c_0^2}=\frac{16}{25},\qquad
 \alpha_1=\frac{j^2}{c_1^2}=\frac{25}{64}.
 \end{gathered}
\]
The triangular form of the block, therefore, gives
\[
 \det(I-D\mathcal Q)=\det(I-A)
 =\prod_{i=1}^k\det\bigl(I_{12}-\alpha_{r_i}W_{K_4}^{\mathsf T}\bigr).
\]

The factors on the right can be evaluated using the determinant polynomial of $K_4$, since transposition does not change a determinant:
\begin{equation}\label{eq:K4poly}
 \det(I-tW_{K_4})
 =(1-t)^3(1+t)^2(1-2t)(1+t+2t^2)^3.
\end{equation}
For completeness, it follows from the ordinary Bass formula and the adjacency eigenvalues $3,-1,-1,-1$: the vertex determinant is $(1-3t+2t^2)(1+t+2t^2)^3$, and the remaining Bass factor is $(1-t^2)^2$. Neither $t=16/25$ nor $t=25/64$ is a root of \eqref{eq:K4poly}. Therefore, $I-D\mathcal Q$ is invertible in each of the $2^k$ solutions.

The analytic implicit-function theorem therefore gives, for each initial solution $\cav_*$ at $\varepsilon=0$, a real-analytic branch $\cav(\varepsilon)$ satisfying $F(\cav(\varepsilon),\varepsilon)=0$ and $\cav(0)=\cav_*$ for $\varepsilon$ near zero. By continuity, the strict inequalities $c_e>0$ and $s_v>0$ persist for sufficiently small $\varepsilon$. The internal output weights start at $4/5$ or $5/8$, while the bridge outputs are positive of the order $\varepsilon$ for positive $\varepsilon$. Thus all output weights remain in $(0,1)$ for sufficiently small positive $\varepsilon$. Since the number of branches is finite, choose one common interval $0<\varepsilon<\varepsilon_k$ on which they exist and remain distinct. Applying \cref{thm:correspondence} completes the proof.
\end{proof}

\begin{remark}\label{rem:manyscope}
The theorem guarantees at least $2^k$ distinct admissible cavity solutions. It does not claim that the constructed solutions are the only ones: other admissible cavity solutions may exist.
The graph has minimum degree at least three, and the directed graph with adjacency matrix $W$ is strongly connected. One may verify this by making a reduced cycle within a $K_4$ before crossing a bridge in either direction. Also $\rho(B(\uu^2))\ge2(4/5)^2>1$, by comparison with the principal submatrix on one copy. Thus the construction satisfies the hypotheses of \cref{thm:main}.
\end{remark}

The construction only requires a finite family of nondegenerate positive cavity solutions on each component, with positive vertex coefficients. By \emph{weak edges} we mean added edges between components, with positive edge weights that tend to zero while all weights inside the original components are held fixed. For sufficiently small values of the added weights, each product choice continues to a cavity solution on the joined graph by the implicit-function theorem. The $K_4$ calculation provides explicit rational data for all hypotheses. It also separates the proof from numerical root finding: the latter is not needed for the asserted existence interval, although no explicit numerical lower bound for $\varepsilon_k$ is claimed here.

\subsection{Example: four branches on two linked copies}\label{sec:fourbranches}
This example illustrates \cref{thm:manybranches} for $k=2$. The graph consists of two labelled copies of $K_4$, on the vertex sets $\{1,\dots,4\}$ and $\{5,\dots,8\}$, joined by the single bridge edge $\{4,5\}$; it consists of the first two blocks of \cref{fig:chain}. Every internal edge carries the weight $4/5$, and the bridge carries the weight $\varepsilon=0.01$. The Newton method in \cref{sec:branchalgorithm}, initialized at the four limiting configurations of \eqref{eq:K4roots}, yields the four distinct cavity solutions summarized in \cref{tab:fourbranches}; the branch label $r_1r_2\in\{00,01,10,11\}$ records which of the constant roots $c_0,c_1$ each copy selects at $\varepsilon=0$.

Each branch determines, through \cref{thm:correspondence}, its own pair $(\uu',S)$ and its own prefactor $C_G$, and for each branch the identity $P_G(\uu')\,C_G=P_G(\uu)$ was checked with the two sides evaluated independently as in \cref{sec:checkalgorithm}. Here $P_G(\uu)=0.206422372892$; the four values of $P_G(\uu')$ are $0.206422372892$, $0.220333813200$ (twice), and $0.235189609178$, and all four residuals are below $10^{-15}$. Branch $00$ is the reference solution, with $\uu'=\uu$ and $C_G=1$. Branch $11$ is the selected branch of \cref{thm:main}, the unique one with $\zz>0$. The equal rows for $01$ and $10$ arise from the symmetry interchanging the copies; they are distinct labelled solutions.
\begin{table}[htbp]
\centering
\caption{Numerical branches at bridge weight $0.01$. The branch label records the root chosen in each copy at $\varepsilon=0$; the near-zero baseline residual is displayed as zero.}\label{tab:fourbranches}
\begin{tabular}{cccccc}
\toprule
Branch & $\min_e z_e$ & $\min_v s_v$ & $\min_e u'_e$ & $\max_e u'_e$ & $C_G$\\
\midrule
$00$ & $0$ & $1$ & $0.01000000$ & $0.80000000$ & $1$\\
$01$ & $-0.00012299$ & $0.99991792$ & $0.00679426$ & $0.80001441$ & $0.93686198$\\
$10$ & $-0.00012299$ & $0.99991792$ & $0.00679426$ & $0.80001441$ & $0.93686198$\\
$11$ & $0.28001595$ & $2.16674552$ & $0.00461552$ & $0.62499003$ & $0.87768492$\\
\bottomrule
\end{tabular}
\end{table}
Each mixed branch has at least one component $z_e<0$, although all its cavity coefficients $c_e$ and vertex coefficients $s_v$ are positive. Thus an admissible solution can lie outside the order region $\cav\ge\cav^{(0)}$. The theorem proves continuation for sufficiently small positive bridge weights; the particular decimal value $\varepsilon=0.01$ is supported by numerical checks, not by an explicitly validated implicit-function radius.

\section{Discussion and further questions}\label{sec:discussion}
The identities in this paper arise from a structured change of edge weights, not from an arbitrary normalization of a zeta quotient. The edge weights and the prefactor are determined by a finite system on directed edges. The cavity parameterization describes all positive diagonal congruences with positive input and output weights less than one, while the nonnegative fixed-point theorem selects and classifies one order region of this parameterization.

Several distinctions delimit the results. First, the nonnegative fixed-point theorem is classical; the SIS correspondence gives an exact route to applying it here. Second, diagonal congruence preserves inertia but not the eigenvalues of $H_G$ or $B(\uu)$. Third, the selected map is defined on $\Omega_+$ and takes values in the complementary squared-subcritical region. It is not presented as a globally rational involution. Finally, all functional identities are first identities of finite determinants. Gaussian integrals are used only when their quadratic forms are positive definite.

The following questions concern structures not classified here. One is the full set of positive cavity solutions with positive vertex coefficients, including their singularities and their connection under changes of the weights. \Cref{thm:manybranches} rules out any universal two-solution answer. A second question is whether the selected map $T_G$ is injective on $\Omega_+$: if $T_G(\uu_1)=T_G(\uu_2)$, the two congruences compose to a positive diagonal congruence between $H_G(\uu_1)$ and $H_G(\uu_2)$, to which \cref{thm:rigidity} does not apply, because neither matrix is positive definite. Another is the subclass of transformations preserving a prescribed family of weights, such as vertex-factorized weights, covering-invariant weights, or a low-dimensional algebraic subvariety of the edge-weight space. Equation \eqref{eq:vertexclosure} isolates the first obstruction but does not classify it globally. A third question is the reducible case of the edge adjacency matrix: its strongly connected components should be retained rather than replacing irreducibility by connectedness of the original graph. Finally, signed, complex, or matrix-valued twists can alter the positivity obstruction; their global branch structure is outside the present scalar-positive theory. A further question is whether $\log C_G$ in \eqref{eq:Cselected} admits an interpretation as a difference of Gaussian Bethe free energies evaluated at the two BP fixed points of \cref{thm:selected}, in the spirit of the orbit-product correction factor of \citet[Theorem~5]{JohnsonChernyakChertkov2009}.

For regular Ramanujan graphs with uniform weights, the graph analogue of the Riemann hypothesis places the nontrivial poles on the critical circle associated with the reciprocal functional equation \cite{Terras2011}. It is natural to ask whether the many-weight correspondence can help formulate analogous constraints for irregular graphs or nonuniform weights. The present results concern positive real weights and preserve the condition $P_G=0$ along the selected transformation; they do not by themselves establish a localization of complex poles. Addressing this question would require control of cavity branches in complex weight space and suitable spectral assumptions. The universal-cover viewpoint of \cref{sec:ratios} may be useful in this direction.

The closest antecedents we have identified are the weighted Ihara--Bass determinant formulas, the Gaussian BP and loop determinant identities, and the theory of ratio systems on universal covers. To the best of our knowledge, the new result is the combination of a parameterization of all positive diagonal congruences within the weight domain $(0,1)^E$ by admissible cavity solutions, the classification in the order region $\cav\ge\cav^{(0)}$ and the resulting real-analytic weight map on $\Omega_+$, rigidity in the positive-definite region, and the construction of connected graphs with exponentially many admissible solutions. \Cref{sec:ratios} gives the precise comparison with ratio-system theory.

\section*{Acknowledgement}
We used ChatGPT (OpenAI, GPT-6) and Claude (Anthropic, Fable 5.1) to assist with drafting and checking the manuscript under our direction. The authors are responsible for all content and conclusions.
The authors would like to thank Kazutoshi Ohta and Ryo Yokokura for useful discussions.
The authors thank the Yukawa Institute for Theoretical Physics at Kyoto University, where this work was initiated during the YITP-W-26-08 on "Strings and Fields 2026".
This work is supported in part by Grant-in-Aid for Scientific Research (KAKENHI) (C), Grant Number 26K07070 (S. M.).

\appendix
\section{Normalized Gaussian kernels and the measure rescaling}\label[appendix]{sec:kernels}
This appendix explains the kernel motivation without making it a prerequisite for the main results.

\subsection{Mehler kernels and the density identity}\label{sec:kerneldensity}
Put
\begin{equation}\label{eq:kernel}
 d\mu(X)=\pi^{-1}e^{-|X|^2}d^2X,\qquad
 K_u(X,Y)=\frac1{1-u^2}
 \exp\left[-\frac{u^2(|X|^2+|Y|^2)}{1-u^2}
 +\frac{u(\overline XY+\overline YX)}{1-u^2}\right].
\end{equation}
Identify $\CC$ with $\RR^2$ and put $u=e^{-t}$, $t>0$. Then $K_u$ is the ground-state-transformed Mehler kernel for the two-dimensional harmonic oscillator \cite{Thangavelu1993}. In the normalization $H_{\mathrm{osc}}=(-\Delta+|X|^2)/2$, the normalized ground state is $\phi_0(X)=\pi^{-1/2}e^{-|X|^2/2}$ and its energy is $1$. The operator $\phi_0^{-1}e^{-t(H_{\mathrm{osc}}-1)}(\phi_0\,\cdot\,)$ has kernel $K_{e^{-t}}$ with respect to $d\mu=\phi_0^2d^2X$. Completing the square gives
\[
 d\mu(Y)K_u(X,Y)
 =\frac{d^2Y}{\pi(1-u^2)}
 \exp\left(-\frac{|Y-uX|^2}{1-u^2}\right),
\]
so $\int K_u(X,Y)d\mu(Y)=1$. Given $X$, this is the law of $Y=e^{-t}X+\sqrt{1-e^{-2t}}\,Z$ with $Z$ distributed according to $\mu$. Thus $K_{e^{-t}}$ is the Ornstein--Uhlenbeck transition density relative to its invariant Gaussian measure $\mu$ \cite{Bogachev2018}. In these variables its semigroup law reads
\[
 \int K_a(X,Y)K_b(Y,Z)d\mu(Y)=K_{ab}(X,Z),\qquad 0<a,b<1.
\]
The scalar Gaussian structure is related to graph-based Kazakov--Migdal representations of graph zeta functions \cite{MatsuuraOhta2022,MatsuuraOhtaWilson2022,MatsuuraOhtaFKM2023,MatsuuraOhta2024}.

Multiplying all edge kernels and all vertex measures yields exactly
\begin{equation}\label{eq:fulldensity}
 \prod_vd\mu(X_v)\prod_{vw\in E}K_{u_{vw}}(X_v,X_w)
 =\frac{\prod_e(1-u_e^2)^{-1}}{\pi^n}
 e^{-X^\dagger H_G(\uu)X}\prod_vd^2X_v.
\end{equation}
Under the change of variables $X'=aX$ with $a>0$, we write the pulled-back Gaussian measure as
\[
 d\mu(aX):=\pi^{-1}e^{-a^2|X|^2}\,a^2\,d^2X.
\]
The factor $a^2$ is the real two-dimensional Jacobian. Applying this change at each vertex with $a=a_v$ and using $S H_G(\uu')S=H_G(\uu)$ gives
\begin{equation}\label{eq:densitytransform}
 \prod_vd\mu(a_vX_v)\prod_{vw\in E}K_{u'_{vw}}(a_vX_v,a_wX_w)
 =C_G\prod_vd\mu(X_v)\prod_{vw\in E}K_{u_{vw}}(X_v,X_w).
\end{equation}
This is a pointwise identity of densities. It can be integrated as an equality of finite positive integrals only in the positive-definite region. \Cref{thm:rigidity} says that, for the positive weights considered in the main text, this region supports only the identity transformation. For a nontrivial selected transformation, the determinant identity is the valid conclusion. When the Gaussian integrals diverge, formally integrating the density identity does not prove the functional equation for the zeta functions.

\subsection{Why positive edge weights are essential for rigidity}\label{sec:signedexample}
The assumption of positive edge weights in \cref{thm:rigidity} is essential. Once negative edge weights are allowed, a nontrivial positive diagonal congruence can exist even when both weighted vertex matrices are positive definite. The following example demonstrates this limitation of the rigidity theorem.

On $K_4$, give every edge the weight $u=-4/5$. Since $K_4$ is $3$-regular, the reciprocal transformation in \eqref{eq:regularrelation}, continued algebraically to negative weights, gives
$u'=\frac{1}{2u}=-\frac58$.
The algebraic weighted Bass formula and the uniform congruence continue to these values, with the corresponding scale $S=\sqrt{13/6}\,I$ and
\[
 S H_G(u'\one)S=H_G(u\one).
\]
Direct diagonalization gives
\[
 \Spec H_G(-4/5)=\{13,37/9,37/9,37/9\},
\]
\[
 \Spec H_G(-5/8)=\{6,74/39,74/39,74/39\}.
\]
Both matrices are positive definite, while $u'\ne u$ and $S\ne I$. Thus positive definiteness alone does not imply rigidity when negative edge weights are admitted. In particular, the signed-weight analogue of \eqref{eq:densitytransform} gives a nontrivial equality of convergent Gaussian integrals. For $0<q<1$ the kernel $K_{-q}$ remains pointwise positive. The identity $K_{-q}(X,Y)=K_q(X,-Y)$ represents an Ornstein--Uhlenbeck step at time $t=-\log q>0$ composed with the reflection $Y\mapsto -Y$. The negative weight $-q$ lies outside the positive-time parametrization $u=e^{-t}$ of the semigroup itself.

The same example also distinguishes convergence of the Gaussian integrals from absolute convergence of the Euler products. Indeed,
\[
 \rho(B(|\uu|))=8/5,\qquad \rho(B(|\uu'|))=5/4,
\]
so the Euler products and cycle expansions do not converge absolutely, although the Gaussian integrals converge. The equivalence in \cref{thm:positivity} therefore also depends on positivity of the edge weights. This example is a signed specialization of the regular-graph reciprocal relation \eqref{eq:regularrelation}; the exact matrix calculation is described in \cref{sec:exampledata}.

\section{Algorithms for reproducing the examples}\label[appendix]{sec:computations}
This appendix specifies the computations independently of a programming language. \Cref{sec:checkalgorithm} gives the common reconstruction and verification procedure; \cref{sec:selectedalgorithm,sec:branchalgorithm} give the two methods for obtaining cavity solutions; and \cref{sec:exampledata} lists the inputs and the exact checks. All arc matrices use \eqref{eq:Wintro}, and every product over $E$ uses each unoriented edge once.

To reproduce the numerical results reported in the main text, use a working precision of $80$ decimal digits. Read finite decimal inputs as exact decimal numbers before converting to the working precision, and round only the reported results. Write $\|x\|_{\max}=\max_i|x_i|$ for a vector and $\|A\|_{\max}=\max_{i,j}|A_{ij}|$ for a matrix. Spectral radii reported numerically are diagnostics computed in double precision. Small residuals in floating-point arithmetic do not certify error bounds or interval enclosures.

\subsection{Common reconstruction and verification}\label{sec:checkalgorithm}
After obtaining a candidate cavity solution by the monotone iteration in \cref{sec:selectedalgorithm} or the Newton iteration in \cref{sec:branchalgorithm}, apply the common procedure below. It constructs the transformed weights, vertex scales, and zeta prefactor from that candidate and checks the cavity equations, the diagonal congruence, and the determinant identities. The two determinant representations, namely, the vertex and edge representations, are evaluated independently, so the verification does not define one side of an identity from the other.

To reproduce the numerical results in the main text, use the verification tolerance $\tau=10^{-58}$.

\begin{algorithm}[Reconstruct and check a cavity solution]\label{alg:checkcavity}
\mbox{}\par\noindent
\textbf{Input:} $G$, reversal-symmetric weights $0<u_e<1$, a candidate $\cav$, and tolerance $\tau$.\par\noindent
\textbf{Output:} transformed data, residuals, and whether the checks pass.
\begin{enumerate}[label=\arabic*.]
\item Construct $W$, $\delta_e=1-u_e^2$, $j_e=u_e/\delta_e$, and $h_v$ from \cref{eq:Wintro,eq:hjintro}. If any $c_e\le0$, report failure.
\item Compute
\[
 s_v=h_v-\sum_{x\sim v}\frac{j_{xv}^2}{c_{x\to v}},\qquad
 u'_{vw}=\frac{j_{vw}}{\sqrt{c_{v\to w}c_{w\to v}}}.
\]
If any $s_v\le0$ or any $u'_e\notin(0,1)$, report failure. Set
\[
 z_e=\delta_ec_e-1,\quad S=\diag(\sqrt{s_v}),\quad
 C=\prod_{v\in V}s_v\prod_{e\in E}\frac{\delta_e}{1-u_e'^2}.
\]
\item Construct $H=H_G(\uu)$ and $H'=H_G(\uu')$ directly from \eqref{eq:HG}. Independently compute
\[
 p=\det(I-WD_{\uu}),\qquad p'=\det(I-WD_{\uu'}),
\]
\[
 b=\prod_{e\in E}\delta_e\det H,\qquad
 b'=\prod_{e\in E}(1-u_e'^2)\det H'.
\]
\item Form the residuals
\[
 \begin{gathered}
 r_{\mathrm{cav}}=\|\cav-\mathcal Q(\cav)\|_{\max},\qquad
 r_{\mathrm{cong}}=\|SH'S-H\|_{\max},\\
 r_{\mathrm{Bass}}=|p-b|,\qquad r'_{\mathrm{Bass}}=|p'-b'|,\qquad
 r_{\mathrm{det}}=|p'-p/C|.
 \end{gathered}
\]
If $p,p'\ne0$, also compute $r_{\zeta}=|1/p'-C/p|$; otherwise omit the reciprocal check.
\item To check the similarity in \cref{prop:image}, set
\[
 J=W^{\mathsf T}D_{(j_e/c_e)^2},\qquad
 L=\diag(s_{o(e)})_{e\in\arcs},\qquad
 r_J=\|L^{-1}JL-W^{\mathsf T}D_{\uu'^2}\|_{\max}.
\]
\item Report whether every computed residual is at most $\tau$. Return $\zz$, $\svec$, $\uu'$, $S$, $C$, $p$, $p'$, and the residuals, together with the diagnostic values $\rho(B(\uu))$, $\rho(B(\uu^2))$, $\rho(B(\uu'))$, and $\rho(B(\uu'^2))$.
\end{enumerate}
\end{algorithm}

\subsection{Monotone iteration for the selected solution}\label{sec:selectedalgorithm}
The same algorithm reproduces
\cref{sec:inhomogeneous} and \cref{ex:diamond}.
It uses \eqref{eq:eta} on all arcs, with no symmetry reduction. Its initial vector is an upper bound because each summand in the update is bounded above by $u_{xv}^2/(1-u_{xv}^2)$. In exact arithmetic the iterates decrease to the selected solution by \cref{lem:iteration,thm:selected}, after the diagonal change between $\eta$ and $\zz$.

To reproduce these examples, use the iteration tolerance $\epsilon_{\mathrm{it}}=10^{-65}$ and the verification tolerance specified in \cref{sec:checkalgorithm}.

\begin{algorithm}[Compute the selected cavity solution]\label{alg:selected}
\mbox{}\par\noindent
\textbf{Input:} $G$ with irreducible $W$, $\uu\in\Omega_+$, $\epsilon_{\mathrm{it}}$, $\tau$, and a maximum of $N_{\mathrm{it}}=20000$ updates.\par\noindent
\textbf{Output:} the selected transformed data and their residuals, or a failure report.
\begin{enumerate}[label=\arabic*.]
\item Set $\delta_e=1-u_e^2$ and initialize
\[
 \eta_{v\to w}=\sum_{x\sim v,\ x\ne w}\frac{u_{xv}^2}{\delta_{xv}}.
\]
\item Repeat at most $N_{\mathrm{it}}$ times. Using the old vector for every component, compute
\[
 \eta^{\mathrm{new}}_{v\to w}
 =\sum_{x\sim v,\ x\ne w}
 \frac{u_{xv}^2\eta_{x\to v}}{1+\delta_{xv}\eta_{x\to v}}.
\]
Set $d=\|\eta^{\mathrm{new}}-\eta\|_{\max}$ and replace $\eta$ by $\eta^{\mathrm{new}}$. If $d<\epsilon_{\mathrm{it}}$, leave the loop. If the update limit is reached without this condition, report failure.
\item Set $c_e=\delta_e^{-1}+\eta_e$. Apply \cref{alg:checkcavity} with tolerance $\tau$, and return its output and the number of updates. Also report whether every computed $z_e$ is positive, as required for the selected branch.
\end{enumerate}
\end{algorithm}

The stopping condition measures the change between successive iterates. The cavity residual is checked separately in \cref{alg:checkcavity}; neither quantity alone supplies a rigorous error enclosure in floating-point arithmetic.

\subsection{Newton iteration for the four cavity branches}\label{sec:branchalgorithm}
To reproduce \cref{sec:fourbranches}, solve at the fixed bridge weight $\varepsilon=0.01$, using the four solutions at $\varepsilon=0$ as initial guesses. This procedure does not require intermediate bridge weights. Define the residual map and its Jacobian by
\begin{equation}\label{eq:Newtonresidual}
 R(\cav):=\cav-\mathcal Q(\cav),\qquad
 DR(\cav)=I-W^{\mathsf T}D_{(j_e/c_e)^2}.
\end{equation}
Here $h_v$, $j_e$, and $\mathcal Q$ always use the target weights, including the positive bridge weight. The limiting roots supply only the starting vectors.

To reproduce this example, use the Newton residual tolerance $\epsilon_{\mathrm N}=10^{-70}$ and the verification tolerance specified in \cref{sec:checkalgorithm}.

\begin{algorithm}[Compute the four branches on two linked copies]\label{alg:branches}
\mbox{}\par\noindent
\textbf{Input:} two copies of $K_4$ on $V_1=\{1,2,3,4\}$ and $V_2=\{5,6,7,8\}$, bridge $\{4,5\}$, internal weight $4/5$, bridge weight $1/100$, $\epsilon_{\mathrm N}$, $\tau$, and a maximum of $N_{\mathrm N}=60$ Newton steps per branch.\par\noindent
\textbf{Output:} checked cavity solutions labelled $00,01,10,11$, or a failure report for a label.
\begin{enumerate}[label=\arabic*.]
\item For each $(r_1,r_2)\in\{0,1\}^2$, perform steps~2--3 with the following initial vector: every internal arc in copy $i$ has value $c_{r_i}$ from \eqref{eq:K4roots}, and the outgoing bridge arc from copy $i$ has value $s_{r_i}$ from \eqref{eq:K4rootdata}.
\item Until $\|R(\cav)\|_{\max}<\epsilon_{\mathrm N}$, perform at most $N_{\mathrm N}$ Newton steps. Solve the linear system
\[
 DR(\cav)\,\Delta=-R(\cav).
\]
Starting with $\lambda=1$, halve $\lambda$ until
\[
 \|R(\cav+\lambda\Delta)\|_{\max}<\|R(\cav)\|_{\max},
\]
and then replace $\cav$ by $\cav+\lambda\Delta$. Report failure for this label if the linear system cannot be solved, the residual cannot be reduced at the working precision, or the step limit is reached without satisfying the residual tolerance.
\item Apply \cref{alg:checkcavity} to the resulting vector and attach the label $r_1r_2$. Record $\min_e z_e$, $\min_v s_v$, $\min_e u'_e$, $\max_e u'_e$, and $C$ for comparison with \cref{tab:fourbranches}. Retain the labelled vectors and check that the four computed solutions are distinct at the working precision.
\end{enumerate}
\end{algorithm}

The mixed branches have negative $z$ components and are not obtained by \cref{alg:selected}. The four successful solves reproduce the displayed branches; they do not enumerate all cavity solutions. Existence for all sufficiently small positive bridge weights is proved by the implicit-function theorem in \cref{thm:manybranches}; the numerical computation at $0.01$ does not give a validated continuation radius.

\subsection{Example data and exact checks}\label{sec:exampledata}
\Cref{tab:algorithmdata} specifies the data for the examples. Decimals in the table are exact input values. Both \cref{alg:selected,alg:branches} include the common checks in \cref{alg:checkcavity}. For the inhomogeneous $K_4$, also compute $u'_{12}u'_{34}$, $u'_{13}u'_{24}$, and $u'_{14}u'_{23}$ from the returned edge weights to reproduce the test of \eqref{eq:K4closure}.

\begin{table}[htbp]
\centering
\caption{Inputs and procedures for reproducing the examples.}\label{tab:algorithmdata}
\begin{tabular}{@{}>{\raggedright\arraybackslash}p{0.23\linewidth}>{\raggedright\arraybackslash}p{0.42\linewidth}>{\raggedright\arraybackslash}p{0.29\linewidth}@{}}
\toprule
Example & Input data & Procedure\\
\midrule
\Cref{ex:k23}: $K_{2,3}$ & $u^2=3/4$, $u'^2=2/3$ & Exact check below\\[4pt]
\Cref{ex:diamond}: diamond & $E=\{12,23,34,41,13\}$; all $u_e=0.85$ & \Cref{alg:selected}\\[4pt]
\Cref{sec:inhomogeneous}: $K_4$ & $q=(0.80,0.85,0.90,0.95)$; $u_{ij}=\sqrt{q_iq_j}$ & \Cref{alg:selected}\\[4pt]
\Cref{sec:fourbranches}: linked $K_4$ copies & Internal weights $4/5$; bridge $\{4,5\}$ of weight $1/100$ & \Cref{alg:branches}\\[4pt]
\Cref{sec:signedexample}: signed $K_4$ & $u=-4/5$, $u'=-5/8$, $S^2=(13/6)I$ & Exact matrix check below\\
\bottomrule
\end{tabular}
\end{table}

The following checks use exact arithmetic; they do not require either numerical root-finding algorithm.
\begin{enumerate}[label=\arabic*.]
\item \textbf{The uniform $K_{2,3}$ example.} Substitute $u^2=3/4$ into \cref{eq:biregularp,eq:biregularscales}, with $r=2$ and $s=1$, to obtain $u'^2=2/3$, $a_A^2=10/7$, and $a_B^2=7/5$. Compute
\[
 C=\left(\frac{10}{7}\right)^2\left(\frac75\right)^3
 \left(\frac{1-3/4}{1-2/3}\right)^6=\frac{5103}{5120}.
\]
Writing $x=u^2$, form $\Phi(x)=(1-x^2)^2(1-4x^2)$ from the uniform Bass formula for $K_{2,3}$. Check $\Phi(2/3)C=\Phi(3/4)$ as an equality of rational numbers.
\item \textbf{The $K_4$ polynomial used in \cref{thm:manybranches}.} Construct the $12\times12$ integer matrix $W_{K_4}$ and compute $D(t)=\det(I-tW_{K_4})$ symbolically. Factor it and compare with \eqref{eq:K4poly}. Evaluate $D(16/25)$ and $D(25/64)$ exactly and check that both are nonzero.
\item \textbf{The signed $K_4$ example.} Construct $H=H_{K_4}(-4/5)$ and $H'=H_{K_4}(-5/8)$ as rational matrices. Check $(13/6)H'=H$ and compute their eigenvalues to reproduce \cref{sec:signedexample}. Using $D$ from the preceding check, set
\[
 C=\left(\frac{13}{6}\right)^4
 \left(\frac{1-16/25}{1-25/64}\right)^6
\]
and verify $D(-5/8)C=D(-4/5)$ exactly. This signed example is checked directly, since the input condition $0<u_e<1$ of \cref{alg:checkcavity} does not apply to it.
\end{enumerate}

\end{document}